\pdfoutput=1
\documentclass[10pt,draft,onecolumn]{ieeetran}
\usepackage[LGR,T1]{fontenc}
\usepackage[latin9]{inputenc}
\usepackage{amsmath}
\usepackage{amssymb,amsfonts}
\usepackage{comment}
\usepackage{epstopdf}

\usepackage{mathtools}
\mathtoolsset{showonlyrefs=true}
\usepackage{color}

\usepackage[noadjust]{cite}
\newtheorem{lemma}{Lemma}
\newtheorem{theorem}{Theorem}
\newtheorem{assumption}{Assumption}

\IEEEoverridecommandlockouts

\title{\LARGE \bf Totally Asynchronous Large-Scale Quadratic Programming: Regularization, Convergence Rates, and Parameter Selection}
\author{Matthew Ubl$^{\star}$ and Matthew T. Hale$^{\star}$
\thanks{$^\star$Department of Mechanical and Aerospace
Engineering, University of Florida, Gainesville, FL, 32611, USA. Emails: \texttt{\{m.ubl,matthewhale\}@ufl.edu}. This research was supported by the office of Naval Research under grant N00014-19-1-2543, the Air Force Office of Scientific Research under the Center of Excellence for Assured Autonomy in Contested Environments, and by a task order from the Munitions Directorate of the Air Force Research Laboratory at Eglin AFB.}}

\begin{document}
\maketitle

\begin{abstract}
Quadratic programs arise in robotics, communications, smart grids,
and many other applications. As these problems grow in size, finding
solutions becomes more computationally demanding, and new algorithms
are needed to efficiently solve them at massive scales. Targeting large-scale problems,
we develop a multi-agent quadratic programming framework in which
each agent updates only a small number of the total decision variables
in a problem. Agents communicate their updated values to each other,
though we do not impose any restrictions on the timing with which
they do so, nor on the delays in these transmissions. Furthermore,
we allow weak parametric coupling among agents, in the sense that
they are free to independently choose their stepsizes, subject to
mild restrictions. We further provide the means
for agents to independently regularize the problems they solve, thereby
improving convergence properties
while preserving agents' independence in selecting parameters and ensuring a global bound on regularization error is satisfied. Larger regularizations accelerate convergence but increase error in the solution obtained, and we quantify the tradeoff between convergence rates and quality of solutions. Simulation results are presented to illustrate these developments.
\end{abstract}

\section{Introduction}
Convex optimization problems arise in a diverse array
of engineering applications, including signal processing~\cite{luo06},
robotics~\cite{schulman14,mellinger11}, 
communications~\cite{chiang05}, 
machine learning~\cite{shalev12},
and many others~\cite{boyd04}. In all of these areas, 
problems can become very large as the number
of network members (robots, processors, etc.) becomes
large.
Accordingly, there has arisen interest in solving
large-scale optimization problems. 
A common feature of large-scale solvers is that they are parallelized
or distributed among a collection of agents in some way. 
As the number of agents grows,
it can be difficult or impossible to ensure synchrony
among distributed computations and communications, and
there has therefore arisen interest in distributed
asynchronous optimization algorithms. 

One line of research considers asynchronous optimization
algorithms in which agents' communication topologies
vary in time. 
A representative sample of this work
includes~\cite{chen12,zhu12,nedic10,nedic09,ram10,lobel11}, 
and these algorithms
all rely on an underlying averaging-based update law, i.e.,
different agents update the same decision variables and then
repeatedly average their iterates to mitigate disagreements
that stem from asynchrony. 
These approaches (and others in the literature)
require some form of graph connectivity over intervals of a
finite length. 
In this paper, we are interested in cases in which delay bounds are
outside agents' control, e.g., 
due to environmental hazards and adversarial jamming for
a team of mobile autonomous agents. 
In these settings,
verifying graph connectivity can be difficult for single
agents to do, and it may not be possible to even check that
connectivity assumptions are satisfied over prescribed intervals. 
Furthermore, even if such checking is possible, it will
be difficult to reliably attain connectivity over the required intervals with
unreliable and impaired communications. For multi-agent systems with
impaired communications, we are interested in developing an
algorithmic framework that succeeds without requiring
delay bounds. 

In this paper, we develop a totally asynchronous quadratic
programming (QP) framework for multi-agent optimization.
Our interest in quadratic programming is motivated
by problems in robotics~\cite{mellinger11} and 
data science~\cite{rodriguez10}, 
where some standard problems can
be formalized as QPs. 
The ``totally asynchronous'' label originates 
in~\cite{bertsekas1989parallel}, and it
describes a class of algorithms which tolerate arbitrarily
long delays, which our framework will do. In addition,
our developments will use block-based update laws in
which each agent updates only a small subset
of the decision variables in a problem, which reduces
each agent's computational burden and, as we will show,
reduces its onboard storage requirements as well. 

Other work
on distributed quadratic programming 
includes~\cite{carli15,teixeira13,lee15,lee2015convergence,kozma2013distributed,todescato2015robust}.
Our work differs from these existing results because we
consider non-separable objective functions, and
because we consider
unstructured update laws (i.e., we do not require communications
and computations to occur in a particular sequence or pattern). 
Furthermore, we consider only deterministic problems, and our framework
converges exactly to a problem's solution, while some existing
works consider stochastic problems and converge approximately or in an appropriate statistical sense. This work is also somewhat related to distributed solutions to systems of linear equations, e.g.~\cite{wang2019solving}, because the gradient of a quadratic function is a linear function. However, methods for such problems are not readily applicable in this paper due to set constraints.

Asynchrony in agents' communications and computations
implies that they will send and receive different information
at different times. As a result, they will disagree about
the values of decision variables in a problem. 
Just as it is difficult for agents to agree on this information,
it can also be difficult to agree on a stepsize value in
their algorithms. One could envision a network of agents solving
an agreement problem, e.g.,~\cite{ren05}, to compute a common 
stepsize,
though we instead allow agents to independently choose
stepsizes, subject to mild restrictions,
thereby eliminating the need to reach agreement
before optimizing. 

It has been shown that
regularizing problems can endow them with an inherent robustness
to asynchrony and improved convergence 
properties, e.g.,~\cite{koshal11,hale15,hale2017asynchronous}. 
Although regularizing is not required here,
we show, in a precise sense, that regularizing improves
convergence rates of our framework as well. 
It is common
for regularization-based approaches to require agents to
use the same regularization parameter, though this is undesirable
for the same reasons as using a common stepsize. Therefore,
we allow agents to independently choose regularization
parameters as well. 

To the best of our knowledge,
few works have considered both independent stepsizes
and regularizations. The most relevant is~\cite{koshal11},
which considers primal-dual algorithms for
problems with functional constraints and synchronous primal updates. 
This paper is different in that we consider
set-constrained problems with totally asynchronous updates, in addition to unconstrained problems.
Regularizing introduces errors in a solution,
and we bound these errors in terms of agents' regularization
parameters.

A preliminary version of this work appeared in~\cite{ubl2019totally}, however this version further includes distributed regularization selection rules for convergence rate and error bound satisfaction, along with new error bounds and and simulation results.

The rest of the paper is organized as follows. 
Section~\ref{sec:background} provides background on QPs
and formal problem statements.
Then, Section~\ref{sec:update} proposes an update law
to solve the problems of interest, and
Section~\ref{sec:convergenceproof} proves its convergence.
Next, Section~\ref{sec:convergencerate} derives a convergence rate, and Section~\ref{sec:regandconvergencerate} then quantifies the effect of regularizations on the convergence rate. Section~\ref{sec:abserrbound} provides an absolute error bound in terms of
agents' regularizations for a set-constrained problem, while Section~\ref{sec:regerrbound} provides a relative error bound for the unconstrained case. Section~\ref{sec:simulation} next illustrates these results in simulation. Finally, Section~\ref{sec:conclusions} concludes the paper.

\section{Background and Problem Statement} \label{sec:background}

In this section, we describe the quadratic optimization
problems to be solved, as well as the assumptions imposed upon these problems
and the agents that solve them. We then describe agents'
stepsizes and regularizations and introduce the need to allow agents
to choose these parameters independently. We next describe the benefits
of independent regularizations, and give two formal problem
statements that will be the focus of the remainder of the paper.
\subsection{Quadratic Programming Background}
We consider a quadratic optimization problem distributed across a
network of $N$ agents, where agents are indexed over $i\in[N]:=\{1,...,N\}$.
Agent $i$ has a decision variable $x^{[i]}\in\mathbb{R}^{n_{i}},n_{i}\in\mathbb{N}$,
which we refer to as its state, and we allow for $n_{i}\neq n_{j}$
if $i\neq j$. The state $x^{[i]}$ is subject to the set constraint
$x^{[i]}\in X_{i}\subset\mathbb{R}^{n_{i}}$, and we make the following
assumption about each $X_{i}$.

\begin{assumption} \label{asm:setconst}
For all $i\in[N]$, the set $X_{i}\subset\mathbb{R}^{n_{i}}$
is non-empty, compact, and convex. $\hfill\triangle$
\end{assumption}

We define the network-level constraint set ${X:=X_{1}\times\cdots\times X_{N}}$,
and Assumption~\ref{asm:setconst} implies that $X$ is non-empty, compact, and convex.
We further define the global state as $x:=\left({x^{[1]}}^{T},...,{x^{[N]}}^{T}\right)^{T}\in X\subset\mathbb{R}^{n}$,
where $n=\sum_{i\in[N]}n_{i}$. We consider quadratic objectives 
\begin{equation}
f(x):=\frac{1}{2}x^{T}Qx+r^{T}x,
\end{equation}
where $Q\in\mathbb{R}^{n\times n}$ and $r\in\mathbb{R}^{n}$. We
then make the following assumption about $f$.
\begin{assumption} \label{asm:Qsymmetric}
In $f$, $Q$ is symmetric. $\hfill\triangle$
\end{assumption}

Note that Assumption~\ref{asm:Qsymmetric} holds without loss of generality because a non-symmetric $Q$ will have only its symmetric part contribute to the value of the quadratic form that defines $f$. Because
$f$ is quadratic, it is twice continuously differentiable, which
we indicate by writing that $f$ is $C^{2}$. In addition, ${\nabla f=Qx+r}$,
and~$\nabla f$ 
is therefore Lipschitz with constant $\|Q\|_{2}$. It is common to assume outright that $Q$ is positive definite, though here we are able to dispense with this assumption based on one in terms of the block structure of agents' updates.

In this paper,
we divide $n\times n$ matrices into blocks. Given a matrix $B\in\mathbb{R}^{n\times n}$,
where $n=\sum_{i=1}^{N}n_{i}$, the $i^{th}$ block of $B$, denoted
$B^{[i]}$, is the $n_{i}\times n$ matrix formed by rows of $B$
with indices $\sum_{k=1}^{i-1}n_{k}+1$ through $\sum_{k=1}^{i}n_{k}$.
In other words, $B^{[1]}$ is the first $n_{1}$ rows of $B$, $B^{[2]}$
is the next $n_{2}$ rows, etc. Similarly, for a vector $b$, $b^{[1]}$
is the first $n_{1}$ entries of $b$, $b^{[2]}$ is the next $n_{2}$
entries, etc. We further define the notation of a sub-block $B^{[i]}_j$, where $B^{[i]} = \left[B^{[i]}_1 \textnormal{ }B^{[i]}_2 \textnormal{ ... } B^{[i]}_N\right]$. That is, $B^{[i]}_{1}$ is the first $n_{1}$ columns of $B^{[i]}$, $B^{[i]}_{2}$ is the next $n_{2}$ columns, etc. For notational simplicity, $B=\left[B^{[i]}_{j}\right]_{p}$ means the matrix $B$ has been partitioned into blocks according to the partition vector $p := [n_{1}, n_{2}, \dots, n_{N}]^{T}$. That is, 
\begin{equation*}
    B = \left[B^{[i]}_{j}\right]_{p}
    = \begin{bmatrix}
    B^{[1]} \\
    B^{[2]} \\
    \vdots \\
    B^{[N]} \\
    \end{bmatrix} =
    \begin{bmatrix}
    B^{[1]}_{1} & B^{[1]}_{2} & \dots & B^{[1]}_{N} \\
    B^{[2]}_{1} & B^{[2]}_{2} & \dots & B^{[2]}_{N} \\
    \vdots & \vdots & \ddots & \vdots \\
    B^{[N]}_{1} & B^{[N]}_{2} & \dots & B^{[N]}_{N} \\
    \end{bmatrix},
\end{equation*}
where $B^{[i]}_{j} \in \mathbb{R}^{n_{i} \times n_{j}}$ for all $i,j \in [N]$

Previous work has shown that totally asynchronous algorithms may diverge if $Q$ is not diagonally dominant~\cite[Example 3.1]{bertsekas1989parallel}. While enforcing this condition is sufficient to ensure a totally asynchronous update scheme will converge, in this paper we will instead require the weaker condition of block diagonal dominance.

\textit{Definition 1:} Let the matrix $B = \left[B^{[i]}_{j}\right]_{p}$, where $p = [n_{1},n_{2},\dots,n_{N}]^{T}$ is given by the dimensions of agents' states above. If the diagonal sub-blocks $B^{[i]}_{i}$ are nonsingular and if 
\begin{equation} \label{eqn:blockdiagdom}
    \left(\left\|B^{[i]^{-1}}_{i}\right\|_{2}\right)^{-1} \geq \sum^{N}_{\substack{j=1 \\ j\neq i}}\left\|B^{[i]}_{j}\right\|_{2} \textnormal{ for all } i \in [N],
\end{equation}
then $B$ is said to be \textit{block diagonally dominant} relative to the choice of partitioning $p$. If strict inequality in Equation~\eqref{eqn:blockdiagdom} is valid for all $i \in [N]$, then $B$ is \textit{strictly block diagonally dominant} relative to the choice of partitioning $p$. \hfill$\blacktriangle$

In later analysis, we will use the gap between the left and right hand side of Equation~\eqref{eqn:blockdiagdom}, which we define as
\begin{equation}
    \delta_{i}(B) = \left(\left\|B^{[i]^{-1}}_{i}\right\|_{2}\right)^{-1} - \sum^{N}_{\substack{j=1 \\ j\neq i}}\left\|B^{[i]}_{j}\right\|_{2}.
\end{equation}

Note that if $p = [1, 1, \dots, 1]^{T}$, Definition 1 reduces to diagonal dominance in the usual sense. We now make the following assumption:
\begin{assumption} \label{asm:Qbdd}
In $f$, $Q=\left[Q^{[i]}_{j}\right]_{p}$ is strictly block diagonally dominant, where $p = [n_{1}, n_{2}, \dots, n_{N}]^{T}$, and $n_{i}$ is the length of $x^{[i]}$ for all $i \in [N]$. $\hfill\triangle$
\end{assumption}

Note also that from Theorem 2 in \cite{feingold1962block}, if Assumptions~\ref{asm:Qsymmetric} and~\ref{asm:Qbdd} hold for a matrix $B$, then $B$ is also positive definite. Therefore Assumptions~\ref{asm:Qsymmetric} and~\ref{asm:Qbdd} imply that $Q \succ 0$, which renders $f$ strongly convex.
\subsection{Problem Statements}

Following our goal of reducing parametric coupling between agents,
we wish to allow agents to select stepsizes independently. In particular, we wish for the stepsize for block $i$, denoted $\gamma_{i}$, to be chosen using only knowledge of $Q^{[i]}$. The selection of $\gamma_{i}$ should not depend on any other block $Q^{[j]}, j \neq i$, or any stepsize choice, $\gamma_{j}$, for any other block.
Allowing independent stepsizes will preclude the need for agents to
agree on a single value before optimizing. The following problem will be one focus of the remainder of
the paper.

\textit{Problem 1:} Design a totally asynchronous distributed optimization
algorithm to solve
\begin{equation}
\underset{x\in X}{\text{minimize}}\quad\frac{1}{2}x^{T}Qx+r^{T}x,
\end{equation}
where only agent $i$ updates $x^{[i]}$, and where agents choose stepsizes
independently. $\hfill\diamond$

While an algorithm that satisfies the conditions stated in Problem
1 is sufficient to find a solution, we wish to allow for regularizations as well.
Regularizations are commonly used for centralized quadratic programs
to improve convergence properties, and we will therefore use
them here. However, in keeping with the independence of agents' parameters,
we wish to allow agents to choose independent regularization parameters. As with stepsizes, we wish for the regularization for block $i$, denoted $\alpha_{i}$, to be chosen using only knowledge of $Q^{[i]}$. The regularized form of $f$, denoted $f_{A}$,
is
\begin{equation} \label{eqn:regform}
f_{A}(x):=f(x)+\frac{1}{2}x^{T}Ax=\frac{1}{2}x^{T}(Q+A)x+r^{T}x,
\end{equation}
where $A=\text{diag}\left(\alpha_{1}I_{n_{1}},...,\alpha_{N}I_{n_{N}}\right)$,
and where $I_{n_{i}}$ is the $n_{i}\times n_{i}$ identity matrix.
Note that $\frac{\partial f_{A}}{\partial x^{[i]}}=Q^{[i]}x+r^{[i]}+\alpha_{i}x^{[i]}$, where 
we see that only $\alpha_{i}$
affects the gradient of $f$ with respect to $x^{[i]}$. With the goal of independent regularizations in mind, we now state
the second problem that we will solve.

\textit{Problem 2:} Design a totally asynchronous distributed optimization
algorithm to solve
\begin{equation}
\underset{x\in X}{\text{minimize}}\quad\frac{1}{2}x^{T}(Q+A)x+r^{T}x,
\end{equation}
where only agent $i$ updates $x^{[i]}$, and where agents independently choose their stepsizes
and regularizations. $\hfill\triangle$

Section III specifies the structure of the asynchronous communications
and computations used to solve Problem 1, and we will solve Problem
1 in Section IV. Afterwards, we will solve Problem 2 in Section V.

\section{Block-Based Multi-Agent Update Law} \label{sec:update}

To define the update law for each agent's state, we first
describe the information stored onboard each agent and how agents
communicate with each other. Each agent will store a vector containing
its own state and that of every agent it communicates with. Formally,
we will denote agent $i$'s full vector of states by $x_{i}$, and
this is agent $i$'s local copy of the global state. Agent $i$'s
own states in this vector are denoted by $x_{i}^{[i]}$. The current
values stored onboard agent $i$ for agent $j$'s states are denoted
by $x_{i}^{[j]}.$ In the forthcoming update law, agent $i$ will only
compute updates for $x_{i}^{[i]}$, and it will share only $x_{i}^{[i]}$
with other agents when communicating. Agent $i$ will only change
the value of $x_{i}^{[j]}$ when agent $j$ sends its own state to agent
$i$. 

At time $k$, agent $i$\textquoteright s full state vector is denoted
$x_{i}(k)$, with its own states denoted $x_{i}^{[i]}(k)$ and those
of agent $j$ denoted $x_{i}^{[j]}(k)$. At any timestep, agent $i$
may or may not update its states due to asynchrony in agents\textquoteright{}
computations. As a result, we will in general have $x_{i}(k)\neq x_{j}(k)$
at all times $k$. We define the set $K^{i}$ to contain all times $k$ at which agent $i$ updates $x_{i}^{[i]}$. In designing
an update law, we must provide robustness to asynchrony while allowing
computations to be performed in a distributed fashion. First-order gradient 
descent methods are robust to many disturbances, with the additional benefit of being computationally simple. Using our notation of a matrix block, we define $\nabla^{[i]}f:=\frac{\partial f}{\partial x^{[i]}}$,
and we see that $\nabla^{[i]}f(x)=Q^{[i]}x+r^{[i]}$, and we propose the following update law: 
\begin{equation}
x_{i}^{[i]}(k+1)=\begin{cases}
\Pi_{X_{i}}\hspace{-0.4em}\left[x_{i}^{[i]}(k)-\hspace{-0.1em}\gamma_{i}\left(Q^{[i]}x_{i}(k)+r^{[i]}\right)\right] & \hspace{-0.5em}k\in K^{i}\\
x_{i}^{[i]}(k) & \hspace{-0.5em}k\notin K^{i}
\end{cases}\hspace{-0.1em},
\end{equation}
where agent $i$ uses some stepsize $\gamma_{i}>0$. The advantage of the block-based 
update law can be seen above, as agent $i$ only needs to know $Q^{[i]}$
and $r^{[i]}$. Requiring each agent to store the entirety of $Q$
and $r$ would require $O(n^{2})$ storage space, while $Q^{[i]}$
and $r^{[i]}$ only require $O(n)$. For large quadratic programs,
this block-based update law dramatically reduces each agent's onboard
storage requirements, which promotes scalability.

In order to account for communication delays, we use $\tau^{j}_{i}(k)$
to denote the time at which the value of $x_{i}^{[j]}(k)$ was originally
computed by agent $j$. For example, if agent $j$ computes a state
update at time $k_{a}$ and immediately transmits it to agent $i$,
then agent $i$ may receive this state update at time $k_{b}>k_{a}$
due to communication delays. Then $\tau^{j}_{i}$ is defined so that
$\tau^{j}_{i}(k_{b})=k_{a}$. For $K^{i}$ and $\tau^{j}_{i}$, we assume
the following.

\begin{assumption} \label{asm:infupdate}
For all $i\in[N]$, the set $K^{i}$ is infinite.
Moreover, for all $i\in[N]$ and $j\in[N]\backslash\{i\}$, if $\left\{ k_{d}\right\} _{d\in\mathbb{N}}$
is a sequence in $K^{i}$ tending to infinity, then
\end{assumption}
\begin{equation}
\lim_{d\rightarrow\infty}\tau^{j}_{i}(k_{d})=\infty. \tag*{$\triangle$}
\end{equation}
Assumption~\ref{asm:infupdate} is simply a formalization of the requirement that no
agent ever permanently stop updating and sharing its own state with
any other agent. For $i\neq j$, the sets $K^{i}$ and $K^{j}$ do
not need to have any relationship because agents' updates are asynchronous.
Our proposed update law for all agents can then be written as follows.

\textit{Algorithm 1:} For all $i\in[N]$ and $j\in[N]\backslash\{i\}$,
execute
\begin{align*}
x_{i}^{[i]}(k+1) & =\begin{cases}
\Pi_{X_{i}}\left[x_{i}^{[i]}(k)-\gamma_{i}\left(Q^{[i]}x_{i}(k)+r^{[i]}\right)\right] & \hspace{-0.5em}k\in K^{i}\\
x_{i}^{[i]}(k) & \hspace{-0.5em}k\notin K^{i}
\end{cases}\\
x_{i}^{[j]}(k+1) & =\begin{cases}
x_{j}^{[j]}\left(\tau^{j}_{i}(k+1)\right) & \hspace{-0.5em}\text{$i$ receives $j$'s state at $k+1$}\\
x_{i}^{[j]}(k) & \hspace{-0.5em}\text{otherwise}\hfill\diamond
\end{cases}
\end{align*}

In Algorithm 1 we see that $x_{i}^{[j]}$ changes only when agent $i$
receives a transmission directly from agent $j$; otherwise it remains
constant. This implies that agent $i$ can update its own state using
an old value of agent $j$\textquoteright s state multiple times and
can reuse different agents\textquoteright{} states different numbers
of times.

\section{Convergence of Asynchronous Optimization}\label{sec:convergenceproof}

In this section, we prove convergence of Algorithm 1. This will be shown using Lyapunov-like convergence.
We will derive stepsize bounds from these concepts that will be used
to show asymptotic convergence of all agents.

\subsection{Block-Maximum Norms}

The convergence of Algorithm 1 will be measured using a block-maximum
norm as 
in~\cite{bertsekas1989convergence},~\cite{bertsekas1989parallel}, 
and~\cite{hale2017asynchronous}. Below, we define the block-maximum
norm in terms of our partitioning vector $p$.

\textit{Definition 2:} Let $x=\left[x^{[i]}\right]_{p}\in\mathbb{R}^{n}$,
where $p = [n_{1},n_{2},\dots,n_{N}]^{T}$. The norm of the full vector $x$ is defined as the maximum 2-norm of
any single block, i.e.,
\begin{equation}
\|x\|_{2,p}:=\max_{i\in[N]}{\|x^{[i]}\|_{2}}. \tag*{$\blacktriangle$}
\end{equation}

The following lemma allows us to upper-bound the induced block-maximum
matrix norm by the norms of the individual blocks. 

\begin{lemma} \label{lem:normsumbound}
For the matrix $B=\left[B^{[i]}_{j}\right]_{p}$ and induced matrix norm $\|B\|_{2,p}$,
\begin{equation}
    \|B\|_{2,p} \leq \max_{i \in [N]}\sum^{N}_{j=1}\left\|B^{[i]}_{j}\right\|_{2}.
\end{equation}
\end{lemma}

\textit{Proof:} Proof in Appendix~\ref{app:normsumbound}. $\hfill\blacksquare$

\subsection{Convergence Via Lyapunov Sub-Level Sets}

We now analyze the convergence of Algorithm 1. We construct
a sequence of sets, $\{X(s)\}_{s\in\mathbb{N}}$, based on work 
in~\cite{bertsekas1989convergence} 
and~\cite{bertsekas1989parallel}. These sets behave analogously to sub-level sets
of a Lyapunov function, and they will enable an invariance type argument
in our convergence proof. Below, we use $\hat{x}:=\arg\min_{x\in X}f(x)$
for the minimizer of $f$. We state the following assumption on these
sets, and below we will construct a sequence of sets that satisfies
this assumption.

\begin{assumption} \label{asm:setsexist}
There exists a collection of sets $\{X(s)\}_{s\in\mathbb{N}}$
that satisfies:

1) $\dots\subset X(s+1)\subset X(s)\subset\dots\subset X$

2) $\lim_{s\rightarrow\infty}X(s)=\left\{ \hat{x}\right\} $

3) There exists $X_{i}(s)\subset X_{i}$ for all $i\in[N]$ and $s\in\mathbb{N}$
such that $X(s)=X_{1}(s)\times...\times X_{N}(s)$

4) $\theta_{i}(y)\in X_{i}(s+1)$, where $\theta_{i}(y):=\Pi_{X_{i}}\left[y^{[i]}-\gamma_{i}\nabla^{[i]}f(y)\right]$
for all $y\in X(s)$ and $i\in[N]$. $\hfill\triangle$
\end{assumption}

Assumptions~\ref{asm:setsexist}.1 and~\ref{asm:setsexist}.2 jointly guarantee that the collection $\{X(s)\}_{s\in\mathbb{N}}$
is nested and that the sets converge to a singleton containing $\hat{x}$.
Assumption~\ref{asm:setsexist}.3 allows for the blocks 
of~$x$ to be updated independently by
the agents, which allows for decoupled update laws. 
Assumption~\ref{asm:setsexist}.4
ensures that state updates make only forward progress toward $\hat{x}$,
which ensures that each set is forward-invariant in time. It is shown
in~\cite{bertsekas1989convergence} and~\cite{bertsekas1989parallel} that the existence of such a sequence of sets
implies asymptotic convergence of the asynchronous update law in Algorithm
1. We therefore use this strategy to show asymptotic convergence
in this paper. We propose to use the construction 
\begin{equation}
X(s)=\left\{ y\in X:\left\|y-\hat{x}\right\|_{2,p}\leq q^{s}D_{o}\right\} ,
\end{equation}
where we define~$D_{o}:=\max_{i\in[N]}\left\|x^{i}(0)-\hat{x}\right\|_{2,p}$,
which is the block furthest from $\hat{x}$ onboard any agent at
timestep zero, and where we define the constant
\begin{equation}
q=\max_{i\in[N]} \left\|I-\gamma_{i}Q^{[i]}_{i}\right\|_{2} + \gamma_{i}\sum^{N}_{\substack{j=1 \\ j\neq i}}\left\|Q^{[i]}_{j}\right\|_{2}.
\end{equation}

To show convergence, we will use the fact that each update contracts towards $\hat{x}$ by a factor
of $q$, and will state a lemma that establishes bounds on every $\gamma_{i}$ 
that imply $q\in(0,1)$. Additionally, we will see that a proof of convergence using this method requires a block diagonal dominance condition on $Q$. This result will be used to show
convergence of Algorithm~1 through satisfaction of 
Assumption~\ref{asm:setsexist}.

If we wish for $q\in(0,1)$, this condition can be restated as
\begin{equation} \label{eqn:igqsum}
    \left\|I-\gamma_{i}Q^{[i]}_{i}\right\|_{2} + \gamma_{i}\sum^{N}_{\substack{j=1 \\ j\neq i}}\left\|Q^{[i]}_{j}\right\|_{2} < 1 \textnormal{ for all } i \in [N].
\end{equation}

Note that because $Q=Q^T \succ 0$ and $Q^{[i]}_{i}$ is a diagonal submatrix of $Q$, we have $Q^{[i]}_{i}=Q^{[i]^{T}}_{i} \succ 0$. From this fact, we see $\left(\left\|Q^{[i]^{-1}}_{i}\right\|_{2}\right)^{-1}=\lambda_{min}\left(Q^{[i]}_{i}\right)$, meaning that Assumption~\ref{asm:Qbdd} holds. Then, in particular,
\begin{equation}
    \lambda_{min}\left(Q^{[i]}_i\right)> \sum^{N}_{\substack{j=1 \\ j\neq i}}\left\|Q^{[i]}_{j}\right\|_{2} \textnormal{ for all } i \in [N].
\end{equation}

The following lemma states an equivalent condition for Equation~\eqref{eqn:igqsum}, which demonstrates the necessity and sufficiency of strict block diagonal dominance.

\begin{lemma} \label{lem:ddnec}
Let $Q=Q^{T}=\left[Q^{[i]}_{j}\right]_{p}$, where $p = [n_{1},n_{2},\dots,n_{N}]^{T}$. Additionally, let the $n \times n$ matrix $\Gamma = \textnormal{ diag}(\gamma_{1}I_{n_{1}},\gamma_{2}I_{n_{2}},...,\gamma_{N}I_{n_{N}})$, where $I_{n_{i}}$ is the identity matrix of size $n_{i}$ and $\gamma_{i}>0$. Then
\begin{equation}
    \left\|I-\gamma_{i}Q^{[i]}_{i}\right\|_{2} + \gamma_{i}\sum^{N}_{\substack{j=1 \\ j\neq i}}\left\|Q^{[i]}_{j}\right\|_{2} < 1 \textnormal{ for all } i \in [N]
\end{equation}
if and only if
\begin{equation}
   \lambda_{min}\left(Q^{[i]}_i\right)> \sum^{N}_{\substack{j=1 \\ j\neq i}}\left\|Q^{[i]}_{j}\right\|_{2}
\textnormal{ and }
    \gamma_{i} < \frac{2}{\sum^{N}_{j=1}\left\|Q^{[i]}_{j}\right\|_{2}}
\end{equation}
for all $i \in [N]$.
\end{lemma}

\textit{Proof:} Proof in Appendix~\ref{app:ddnec}. $\hfill\blacksquare$

Note that $\gamma_{i}$ only depends on $Q^{[i]}$. This lemma implies that $\gamma_{i}$ can be chosen according to the conditions of Problem 1 such that $q\in(0,1)$, given that Assumption~\ref{asm:Qbdd} holds for $Q$. Choosing appropriate stepsizes for all $i\in[N]$ and recalling our construction of sets $\left\{ X(s)\right\} _{s\in\mathbb{N}}$
as
\begin{equation} \label{eqn:setcon}
X(s)=\left\{ y\in X:\|y-\hat{x}\|_{2,p}\leq q^{s}D_{o}\right\} ,
\end{equation}
we next show that Assumption~\ref{asm:setsexist} is satisfied, thereby ensuring
convergence of Algorithm 1.

\begin{theorem} \label{thm:setswork}
If Assumptions 1-4 hold and $\Gamma = \textnormal{ diag}(\gamma_{1}I_{n_{1}},\gamma_{2}I_{n_{2}},...,\gamma_{N}I_{n_{N}})$ satisfies the conditions in Lemma~\ref{lem:ddnec}, then the collection of sets $\left\{ X(s)\right\} _{s\in\mathbb{N}}$
as defined in Equation~\eqref{eqn:setcon} satisfies Assumption~\ref{asm:setsexist}.
\end{theorem}
\textit{Proof:} Proof in Appendix~\ref{app:setswork}. $\hfill\blacksquare$

Regarding Problem 1, we therefore state the following:

\begin{theorem} \label{thm:alg1works}
Algorithm 1 solves Problem 1 and asymptotically
converges to $\hat{x}$.
\end{theorem}

\textit{Proof:} Proof in Appendix~\ref{app:alg1works}. $\hfill\blacksquare$

From these requirements, we see that agent $i$ only needs to be
initialized with $Q^{[i]}$ and $r^{[i]}$. Agents are then
free to choose stepsizes independently, provided
stepsizes obey the bounds established in Lemma~\ref{lem:ddnec}. 

\section{Convergence Rate} \label{sec:convergencerate}

Beyond asymptotic convergence, the structure of the sets $\left\{ X(s)\right\} _{s\in\mathbb{N}}$ allows us to determine a convergence rate. To do so, we first define the notion of a \textit{communication cycle}. 

\textit{Definition 3:} One \textit{communication cycle} occurs when every agent has calculated a state update and this updated state has been sent to and received by every other agent.\hfill$\blacktriangle$

Once the last updated state has been received by the last agent, a communication cycle ends and another begins. It is only at the conclusion of the first communication cycle that each agents' copy of the ensemble state is moved from $X(0)$ to $X(1)$. Once another cycle is completed every agent's copy of the ensemble state is moved from $X(1)$ to $X(2)$. This process repeats indefinitely, and coupled with Assumption~\ref{asm:setsexist}, means the convergence rate is geometric in the number of cycles completed, which we now show.

\begin{theorem} \label{convergencerate}
Let Assumptions 1-5 hold and let $\gamma_{i} \in \left(0,\frac{2}{\sum^{N}_{j=1}\left\|Q^{[i]}_{j}\right\|_{2}}\right)$ for all $i \in [N]$. At time $k$, if $c(k)$ cycles have been completed, then $\|x_{i}(k)-\hat{x}\|_{2,p} \leq q^{c(k)}D_{o}$
for all $i \in [N]$.
\end{theorem}

\textit{Proof:} Proof in Appendix~\ref{app:convergencerate}. $\hfill\blacksquare$

From the definition of $q$, we may write
    $q = \max_{i\in[N]}q_{i}$,
where
\begin{equation} \label{eqn:defqi}
    q_{i} =  \left\|I-\gamma_{i}Q^{[i]}_{i}\right\|_{2} + \gamma_{i}\sum^{N}_{\substack{j=1 \\ j\neq i}}\left\|Q^{[i]}_{j}\right\|_{2},
\end{equation}
which illustrates the dependence of each $q_{i}$ upon $\gamma_{i}$. As in all forms of gradient descent optimization, the choice of stepsizes has a significant impact on the convergence rate, which can be expressed through its effect on $q$. Therefore, we would like to determine the optimal stepsizes for each block in order to minimize $q$, which will accelerate convergence to a solution. Due to the structure of $q$, minimizing $q_{i}$ for each $i\in[N]$ will minimize $q$. This fact leads to the following theorem:
\begin{theorem} \label{thm:optimalstep}
$q$ is minimized when, for every $i\in[N]$, 
\begin{equation}
     \gamma_{i} = \frac{2}{\lambda_{max}\left(Q^{[i]}_i\right)+\lambda_{min}\left(Q^{[i]}_i\right)}.
\end{equation}
\end{theorem}

\textit{Proof:} Proof in Appendix~\ref{app:optimalstep}. $\hfill\blacksquare$

\section{Regularization and Convergence Rate} \label{sec:regandconvergencerate}
In centralized optimization, regularization can be used to accelerate convergence by reducing the condition number of $Q$. It is well known that the condition number of $Q$, denoted $k_{Q}$, plays a significant role in the convergence rate, with large condition numbers correlating to slow convergence rates. However in a decentralized setting it is difficult for agents to independently select regularizations such that $k_Q$ is reduced, and harder still to know the magnitude of the reduction. In~\cite{ubl2019totally} it is shown that if the ratio of the largest to smallest regularization used in the network is less than $k_Q$, then the condition number of the regularized problem is guaranteed to be smaller. However, this requires global knowledge of $k_Q$, requires an upper bound on regularizations to somehow be agreed on, and institutes a lower bound on agents' choice of regularizations, all of which lead to the type of parametric coupling that we wish to avoid. 

As stated in
Problem 2, we want to allow agents to choose regularization parameters
independently. Here, we therefore only require that agent $i$ use a positive regularization parameter $\alpha_{i} > 0$. In Algorithm 1, this changes only agent $i$'s updates to $x^{[i]}_{i}$, which now take the form
\begin{equation}
    x^{[i]}_{i} = \Pi_{X_{i}}\left[x_{i}^{[i]}(k)-\gamma_{i}\left(Q^{[i]}x_{i}(k)+r^{[i]}+\alpha_{i}x^{[i]}_{i}(k)\right)\right].
\end{equation}
Before we analyze the effects of independently chosen
regularizations on convergence, we must first show that an algorithm
that utilizes them will preserve the convergence properties of Algorithm~1. As shown in Equation~\eqref{eqn:regform}, a regularized cost function takes the form
\begin{equation}
f_{A}(x):=\frac{1}{2}x^{T}(Q+A)x+r^{T}x,
\end{equation}
where $Q+A$ is symmetric and positive definite because ${Q=Q^{T}\succ0}$.
We now state the following theorem that confirms that minimizing $f_{A}$
succeeds.

\begin{theorem} \label{thm:prob2solved}
Suppose that $A=\text{diag}\left(\alpha_{1}I_{n_{1}},...,\alpha_{N}I_{n_{N}}\right)\succ0$, where agent $i$ chooses $\alpha_{i}$ independently of all other agents. 
Then Algorithm
1 satisfies the conditions stated in Problem 2 when $f_{A}$ is minimized.
\end{theorem}

\textit{Proof:} Replacing $Q$ with $Q+A$, all assumptions and conditions used to prove Theorem~\ref{thm:alg1works} hold, with the only modifications being the network will converge to $\hat{x}_{A} := \arg\min_{x\in X} f_{A}(x)$. These steps are similar to those used to prove Theorem~\ref{thm:alg1works} and are therefore omitted. \hfill$\blacksquare$

Theorem~\ref{thm:prob2solved} establishes that regularizing preserves asymptotic convergence, and we next turn to analyzing convergence rates. Because the condition number $k_{Q}$ is a parameter that depends on the entirety of $Q$, and each agent only has access to a portion of $Q$, it is impossible for agents to know how their independent choices of regularizations affect $k_{Q}$. However, we can instead use $q$, which provides our convergence rate and can be directly manipulated by agents' choice of regularizations. Assume the optimal stepsize for block $i$ is chosen as given in Equation~\eqref{eqn:optstep}. We then have
\begin{equation}
    q_{i}  = \frac{2\sum^{N}_{j\neq i}\left\|Q^{[i]}_{j}\right\|_{2}+ \lambda_{max}\left(Q^{[i]}_i\right)-\lambda_{min}\left(Q^{[i]}_i\right)}{\lambda_{max}\left(Q^{[i]}_i\right)+\lambda_{min}\left(Q^{[i]}_i\right)}.
\end{equation}

When we regularize the problem with $A$, the convergence parameter becomes $q_{A} = \max_{i}q_{\alpha_{i}}$, where
\begin{align*}
    q_{\alpha_i}
    & = \frac{2\sum^{N}_{j\neq i}\left\|Q^{[i]}_{j}\right\|_{2}+ \lambda_{max}\left(Q^{[i]}_i\right)-\lambda_{min}\left(Q^{[i]}_i\right)}{\lambda_{max}\left(Q^{[i]}_i\right)+\lambda_{min}\left(Q^{[i]}_i\right)+2\alpha_{i}}.
\end{align*}

The only effect regularization has on $q_{i}$ is adding $2\alpha_i$ to the denominator, meaning that \textit{any} choice of positive regularization will result in $q_{\alpha_{i}} < q_i$, and thus all regularizations accelerate convergence. Using this fact, we can tailor parameter selections to attain a desired convergence rate. Assume we have a desired convergence rate for our system, corresponding to $q^*$. If we want to set $q_{A} \leq q^{*}$, we need $q_{\alpha_{i}} \leq q^{*}$ for all $i \in [N]$. Some algebraic manipulation of the above equation shows we therefore need to choose $\alpha_{i}$ such that
\begin{equation}
    \alpha_{i} \geq \left(\frac{q_{i}}{q^{*}}-1\right)\left(\frac{\lambda_{max}\left(Q^{[i]}_{i}\right)+\lambda_{min}\left(Q^{[i]}_{i}\right)}{2}\right).
\end{equation}

Note that this term will be negative if $q_{i}<q^{*}$. That is, if the dynamics of block $i$ are such that it will already converge faster than required by $q^{*}$, then there is no need to regularize that block.
We now state the following theorem:

\begin{theorem} \label{thm:qreg}
Given $q^{*} \in (0,1)$, if for all $i\in[N]$ agent $i$ chooses
\begin{equation} \label{eqn:qgammabound}
    \gamma_{i} = \frac{2}{\lambda_{max}\left(Q^{[i]}_i\right)+\lambda_{min}\left(Q^{[i]}_i\right)+2\alpha_{i}},
\end{equation}
where
\begin{equation}
     \alpha_{i} = \max\left\{\hspace{-0.3em}\left(\frac{q_{i}}{q^{*}}-1\right)\hspace{-0.3em}\left(\frac{\lambda_{max}\left(Q^{[i]}_{i}\right)+\lambda_{min}\left(Q^{[i]}_{i}\right)}{2}\right),0\hspace{-0.1em}\right\},
\end{equation}
then $q_{A}\leq q^{*}$.
\end{theorem}
\textit{Proof:} Substitute Equation~\eqref{eqn:qgammabound} into Equation~\eqref{eqn:defqi}. \hfill$\blacksquare$

\section{Regularization Absolute Error Bound: Set Constrained Case}
\label{sec:abserrbound}
Regularization inherently results in a suboptimal solution because the system converges to $\Pi_{X}\left[\hat{x}_{A}\right]$ rather than $\Pi_{X}\left[\hat{x}\right]$. We therefore wish to bound this error by a function of the regularization matrix $A$. We define this error in two ways, $\left\|\Pi_{X}\left[\hat{x}\right]-\Pi_{X}\left[\hat{x}_{A}\right]\right\|_{2,p} = \max_{i}\left\|\Pi_{X_{i}}\left[\hat{x}^{[i]}\right]-\Pi_{X_{i}}\left[\hat{x}_{A}^{[i]}\right]\right\|_{2}$, which is the largest error of any one block in the network, and $\left|f\left(\Pi_{X}\left[\hat{x}\right]\right)-f\left(\Pi_{X}\left[\hat{x}_{A}\right]\right)\right|$, which is the difference in cost for the system between the regularized and unregularized cases. Note that in this section we are deriving descriptive error bounds in the sense that a network operator with access to each agent's local information can bound the error for the entire system, but no individual agent is expected to have access to this information.

Looking at the first definition of error, we find
\begin{equation}
    \left\|\Pi_{X}\left[\hat{x}\right]-\Pi_{X}\left[\hat{x}_{A}\right]\right\|_{2,p} \leq \left\|\hat{x}-\hat{x}_{A}\right\|_{2,p}
\end{equation}
which follows from the non-expansive property of the projection operator. Because of the fact that $\hat{x} = -Q^{-1}r$ and ${\hat{x}_{A} = -(Q+A)^{-1}r}$, we see
\begin{equation}
    \left\|\hat{x}-\hat{x}_{A}\right\|_{2,p} = \left\|(Q^{-1}-(Q+A)^{-1})r\right\|_{2,p}.
\end{equation}
Through use of the Woodbury matrix identity, one can see $Q^{-1}-(Q+A)^{-1} = (I+A^{-1}Q)^{-1}Q^{-1}$, because $A$ is invertible. This gives 
\begin{align} \label{eqn:x-xaform}
    \left\|\hat{x}-\hat{x}_{A}\right\|_{2,p}
    & \leq \left\|(I+A^{-1}Q)^{-1}\right\|_{2,p}\left\|Q^{-1}\right\|_{2,p}\left\|r\right\|_{2,p}.
\end{align}
Here $\left\|r\right\|_{2,p} = \max_{i}\left\|r^{[i]}\right\|_{2}$ is the largest norm of any individual block of $r$, which a network operator can gather from agents. However, the two other terms are $2,p$-norms of inverse matrices, which we do not assume the network operator has the ability to calculate. However, these terms can be bounded above using local information from agents according to the following lemma.
\begin{lemma} \label{lem:invinftynorm}
If there is a block strictly diagonally dominant matrix $B = \left[B^{[i]}_{j}\right]_{p}$, where $p = [n_{1},n_{2},\dots,n_{N}]^{T}$, and $\beta_{p}(B) = \min_{i}\left(\left\|B^{[i]^{-1}}_{i}\right\|_{2}^{-1}-\sum^{N}_{\substack{j=1 \\ j\neq i}}\left\|B^{[i]}_{j}\right\|_{2}\right)$, then
\begin{equation}
    \left\|B^{-1}\right\|_{2,p} \leq \beta^{-1}_{p}(B).
\end{equation}
\end{lemma}
\textit{Proof:} Theorem 2 in~\cite{varah1975lower} establishes the above result for $\|\cdot\|_{\infty}$, and the proof for $\|\cdot\|_{2,p}$ follows identical steps. $\hfill\blacksquare$

We note also that $I+A^{-1}Q$ is strictly block diagonally dominant, as $(A^{-1}Q)^{[i]} = \alpha_{i}^{-1}Q^{[i]}$. That is, each block of $Q$ is multiplied by a positive scalar, which preserves the strict diagonal dominance of each block, as does the addition of $I$. Therefore, using Lemma~\ref{lem:invinftynorm} and $Q^{[i]}_{i}=Q^{[i]^{T}}_{i} \succ 0$ for all $i \in [N]$ we see $\left\|(I+A^{-1}Q)^{-1}\right\|_{2,p} \leq \beta^{-1}_{p}(I+A^{-1}Q)$
and $\left\|Q^{-1}\right\|_{2,p} \leq \beta^{-1}_{p}(Q)$,
where $\beta_{p}(I+A^{-1}Q) = \min_{i}\left(1+\alpha^{-1}_{i}
    \delta_{i}(Q)
    \right)$
and 
$\beta_{p}(Q) = \min_{i}\delta_{i}(Q).$ Finally,
\begin{equation} \label{eqn:absxbound}
    \left\|\Pi_{X}\left[\hat{x}\right]-\Pi_{X}\left[\hat{x}_{A}\right]\right\|_{2,p} \leq \frac{\max_{i}\|r^{[i]}\|_{2}}{\beta_{p}(I+A^{-1}Q)\beta_{p}(Q)}.
\end{equation}
The significance of this error bound is that if a network operator has access to $\|r^{[i]}\|_{2}$,  $\alpha_{i}$, and $\delta_{i}(Q)$ for all $i \in [N]$, which are locally known to every agent, the network operator can compute these bounds.

Defining $\Delta_{X_{A}} = \Pi_{X}\left[\hat{x}\right]-\Pi_{X}\left[\hat{x}_{A}\right]$, we find that $f(\Pi_{X}\left[\hat{x}\right])-f(\Pi_{X}\left[\hat{x}_{A}\right]) = \frac{1}{2}(\Pi_{X}\left[\hat{x}\right]+\Pi_{X}\left[\hat{x}_{A}\right])^{T}Q(\Delta_{X_{A}})+r^{T}(\Delta_{X_{A}})$, which gives
\begin{align}
    &|f(\Pi_{X}\left[\hat{x}\right])-f(\Pi_{X}\left[\hat{x_{A}}\right])| \\
    & = \Big|\frac{1}{2}(\Pi_{X}\left[\hat{x}\right]+\Pi_{X}\left[\hat{x}_{A}\right])^{T}Q(\Delta_{X_{A}}) +r^{T}(\Delta_{X_{A}})\Big| \\
    & \leq \|\frac{1}{2}(\Pi_{X}\left[\hat{x}\right]+\Pi_{X}\left[\hat{x}_{A}\right])^{T}Q+r^{T}\|_{2,p}\|\Delta_{X_{A}}\|_{2,p} \\
    & \leq (\|\frac{1}{2}(\Pi_{X}\left[\hat{x}\right]+\Pi_{X}\left[\hat{x}_{A}\right])^{T}Q\|_{2,p}+\|r^{T}\|_{2,p})\|\Delta_{X_{A}}\|_{2,p} \\
    & \leq (\|\frac{1}{2}(\Pi_{X}\left[\hat{x}\right]\hspace{-0.1em}+\hspace{-0.1em}\Pi_{X}\left[\hat{x}_{A}\right])^{T}\|_{2,p}\|Q\|_{2,p}\hspace{-0.1em}+\hspace{-0.1em}\|r^{T}\|_{2,p})\|\Delta_{X_{A}}\|_{2,p}.
\end{align}
Note that by definition, $\|x^{T}\|_{2,p} = \sum^{N}_{i=1}\|x^{[i]}\|_{2}$, and by Lemma 1 $\|B\|_{2,p} \leq \max_{i}\sum^{N}_{j=1}\left\|B^{[i]}_{j}\right\|_{2}$. Combining this with the non-expansive property of the projection operator gives
\begin{align}
    |f(\Pi_{X}\left[\hat{x}\right])-f(\Pi_{X}\left[\hat{x_{A}}\right])| \\
    \leq \bigg(\max_{i}\Big\|\frac{1}{2}\Big(\Pi_{X_{i}}\left[\hat{x}^{[i]}\right]+\Pi_{X_{i}}&\left[\hat{x}^{[i]}_{A}\right]\Big)^{T}\Big\|_{2}\max_{i}\sum^{N}_{j=1}\left\|Q^{[i]}_{j}\right\|_{2} \\
    & +\sum^{N}_{i=1}\left\|r^{[i]}\right\|_{2}\bigg)\|\Delta_{X_{A}}\|_{2,p}.
\end{align}
From Assumption 1, the set constraint for each block is compact, meaning agents can find the vector $\bar{x}^{[i]} = \textnormal{arg}\max_{x^{[i]} \in X_{i}} \|x^{[i]}\|_{2}$. Setting $M_{X_{i}} = \|\bar{x}^{[i]}\|_{2}$ and combining this with Equation~\eqref{eqn:absxbound} gives
\begin{align}
    &|f(\Pi_{X}\left[\hat{x}\right])-f(\Pi_{X}\left[\hat{x_{A}}\right])| \\
    &
    \leq \frac{(\max_{i \in [N]}M_{X_{i}}\max_{i \in [N]}\sum^{N}_{j=1}\left\|Q^{[i]}_{j}\right\|_{2}+\sum^{N}_{i=1}\|r^{[i]}\|_{2})}{\beta_{p}(I+A^{-1}Q)\beta_{p}(Q)} \\
    & + \frac{\max_{i \in [N]}\|r^{[i]}\|_{2}}{\beta_{p}(I+A^{-1}Q)\beta_{p}(Q)}.
\end{align}

\section{Regularization Relative Error Bound: Unconstrained Case} \label{sec:regerrbound}

In the previous section we derived a descriptive bound for the absolute error in both the states of the system and the cost due to regularizing. This bound is descriptive in the sense that given the agents' regularization choices, one can derive a bound describing error for the system. However given a desired error bound, agents cannot use the above rules to independently select regularizations due to the need for global information. Eliminating this dependence upon global information appears to be difficult because of the wide range of possibilities for the set constraints $X_{i}$. However, in the case where our problem does not have set constraints, i.e. Assumption 1 no longer holds and $X = \mathbb{R}^{n}$, we find that we can develop an entirely independent regularization selection rule to bound relative error. In particular, given some $\epsilon > 0$, we wish to bound the relative cost error via
\begin{equation}
    \frac{\left| f(\hat{x}) - f(\hat{x}_{A})\right|}{\left| f(\hat{x})\right|} \leq \epsilon.
\end{equation}

If agents independently select regularizations, then $\alpha_{i}$ is selected using only knowledge of $Q^{[i]}$. Because we do not want to require agents to coordinate their regularizations to ensure the error bound is satisfied, we must develop independent regularization selection guidelines that depend only on $Q^{[i]}$.

\textit{Problem 3:} Given the restriction that $\alpha_{i}$ can be chosen using only knowledge of $Q^{[i]}$ and $\epsilon$, where $\epsilon \in (0,1)$, develop independent regularization selection guidelines that guarantee
\begin{equation}
    \frac{\left| f(\hat{x}) - f(\hat{x}_{A})\right|}{\left| f(\hat{x})\right|} \leq \epsilon. \tag*{$\triangle$}
\end{equation}

For the unregularized problem, the solution is $\hat{x} = -Q^{-1}r$ and the optimal cost is $f(\hat{x}) = -\frac{1}{2}r^{T}Q^{-1}r$. For the regularized problem, the regularized solution is $\hat{x}_{A} = -P^{-1}r$, where $P = Q+A$, and the suboptimal cost is $f(\hat{x}_{A}) = \frac{1}{2}r^{T}P^{-1} QP^{-1}r - r^{T}P^{-1}r$. Note that $f(\hat{x}) \leq f(\hat{x}_{A}) \leq 0$. That is, the cost can be upper-bounded by zero trivially for both the regularized and unregularized cases using $x=0$. Therefore the optimal cost in both cases will be negative, with $f(\hat{x}) \leq f(\hat{x}_{A})$. In particular, we know $f(\hat{x}) - f(\hat{x}_{A}) \leq 0$ and $f(\hat{x}) \leq 0$. Assuming $f(\hat{x}) \neq 0$, we can say 
\begin{equation}
     \frac{ f(\hat{x}) - f(\hat{x}_{A})}{f(\hat{x})} \geq 0.
\end{equation}

That is,
\begin{equation}
    \frac{\left| f(\hat{x}) - f(\hat{x}_{A})\right|}{\left| f(\hat{x})\right|} \leq \epsilon \textnormal{  if and only if  } \frac{ f(\hat{x}) - f(\hat{x}_{A})}{f(\hat{x})} \leq \epsilon.
\end{equation}

The solution to Problem 3 will be developed in two parts. First, it will be shown that the block diagonal dominance condition of $Q$ allows $A$ to be chosen under the restrictions of Problem 3 such that a certain eigenvalue condition of the matrix $A^{-1}Q$ is satisfied. Afterward, it will be shown that this condition on $A^{-1}Q$ is sufficient to guarantee the error bound given by $\epsilon$ is satisfied.

\subsection{Block Gershgorin Circle Theorem}

The Gershgorin Circle Theorem tells us that for any eigenvalue of a symmetric $n \times n$ matrix $B$, we have ${\lambda_{k}(B) \in \bigcup_{k=1}^{n}[b_{k,k} - \sum^{n}_{j \neq k}|b_{k,j}|,b_{k,k} + \sum^{n}_{j \neq k}|b_{k,j}|]}$ for all ${k = 1,...,n}$. That is, every eigenvalue of $B$ is contained within a union of intervals dependent on the rows of $B$. This implies that we can lower bound the minimum eigenvalue of $B$ by ${\lambda_{min}(B) \geq \min_{k}(b_{k,k} - \sum^{n}_{j \neq k}|b_{k,j}|)} $. In the event that $B$ is a strictly diagonally dominant matrix in the usual sense, i.e., $n_{i} = 1$ for all $i \in [N]$, this implies that every eigenvalue of $B$ is positive, because ${\lambda_{min}(B) \geq \min_{k}b_{k,k} - \sum^{n}_{j \neq k}|b_{k,j}|} > 0$ for all ${k = 1,...,n}$. Note further that if we let $C$ be an $n \times n$ positive definite diagonal matrix, then ${\lambda_{min}(CB) \geq \min_{k}c_{k,k}(b_{k,k} - \sum^{n}_{j \neq k}|b_{k,j}|) > 0}$. That is, if $B$ is a strictly diagonally dominant matrix and $C$ is a positive definite diagonal matrix, then $CB$ is strictly diagonally dominant. 

Let $B$ and $C$ meet the criteria above, and now let us treat $C$ as a design choice. Suppose we wish for the smallest eigenvalue of $CB$ to be greater than or equal to a particular constant $l$, i.e., we want $\lambda_{min}(CB) \geq l$. From the Gershgorin Circle Theorem, we see this is true if $c_{k,k}(b_{k,k} - \sum^{n}_{j \neq k}|b_{k,j}|) \geq l$ for all ${k = 1,...,n}$. This condition can be restated as
\begin{equation}
    \textnormal{if } c_{k,k} \geq \frac{l}{b_{k,k} - \sum^{n}_{j \neq k}|b_{k,j}|} \textnormal{ for all } k = 1,...,n,
\end{equation}
then $\lambda_{min}(CB) \geq l$.

That is, given a strictly diagonally dominant matrix $B$ and a positive constant $l$, the $k^{th}$ diagonal element of $C$ can be chosen using only knowledge of the $k^{th}$ row of $B$ and $l$ such that $\lambda_{min}(CB) \geq l$. This intuition can be extended to a strictly block diagonally dominant matrix $B$ using a block analogue of the Gershgorin Circle Theorem, as described below.

\begin{lemma} \label{lem:blockgersh}
For the matrix $B = \left[B^{[i]}_{j}\right]_{p}$, where $p = [n_{1},n_{2},\dots,n_{N}]^{T}$, each eigenvalue $\lambda(B)$ satisfies
\begin{equation}
    \left(\left\|\left(B^{[i]}_{i}-\lambda(B) I\right)^{-1}\right\|_{2}\right)^{-1} \leq \sum^{N}_{\substack{j=1 \\ j\neq i}}\left\|B^{[i]}_{j}\right\|_{2}
\end{equation}
for at least one $i \in [N]$.

\end{lemma}

\textit{Proof:} See Theorem 2 in~\cite{feingold1962block}. $\hfill\blacksquare$

Note that 
\begin{equation}
\left(\left\|\left(B^{[i]}_{i}\hspace{-0.1em}-\hspace{-0.1em}\lambda_{min}(B) I\right)^{-1}\right\|_{2}\right)^{-1} \hspace{-0.75em}= \min_{i}\left|\lambda_{min}(B)\hspace{-0.1em}-\hspace{-0.1em}\lambda_{i}\left(B^{[i]}_{i}\right)\right|\hspace{-0.1em}.
\end{equation}
Additionally, let 
\begin{equation}
    \mu\left(B^{[i]}_{i}\right) = \arg\min_{\lambda_{i}}\left|\lambda_{min}(B)-\lambda_{i}\left(B^{[i]}_{i}\right)\right|,
\end{equation} which is the eigenvalue of $B^{[i]}_{i}$ closest to the minimum eigenvalue of $B$. Then,
\begin{equation}
\left(\left\|\left(B^{[i]}_{i}-\lambda_{min}(B) I\right)^{-1}\right\|_{2}\right)^{-1} = \left|\lambda_{min}(B)-\mu\left(B^{[i]}_{i}\right)\right|. 
\end{equation}

From the block Gershgorin Circle Theorem, we then have
\begin{equation}
    \lambda_{min}(B) \geq \mu\left(B^{[i]}_{i}\right) - \sum^{N}_{\substack{j=1 \\ j\neq i}}\left\|B^{[i]}_{j}\right\|_{2} \textnormal{for at least one } i \in [N].
\end{equation}

Because $\mu\left(B^{[i]}_{i}\right) \geq \lambda_{min}\left(B^{[i]}_{i}\right)$, we can say
$\lambda_{min}(B) \geq \delta_{i}(B) \textnormal{ for at least one } i \in [N].$

Just as before, if $B$ is strictly block diagonally dominant, then every eigenvalue of $B$ is positive. Now let $C=\left[C^{[i]}_{j}\right]_{p}$, with $C^{[i]}_{i} = c_{i}I$ for every $i \in [N]$ and $C^{[i]}_{j} = 0$ when $j \neq i$. In the same manner as above, we find
\begin{equation} \label{eqn:CB>l}
    \textnormal{if } c_{i} \geq \frac{l}{\delta_{i}(B)} \textnormal{ for all } i \in [N],
\end{equation}
then $\lambda_{min}(CB) \geq l$.

That is, $c_{i}$ can be chosen using only knowledge of $B^{[i]}$ and $l$. This brings us back to the restrictions imposed in Problem 3. For reasons that will be shown in the following subsection, choose $B=Q$, $C = A^{-1}$, and $l = \frac{1-\sqrt{\epsilon}}{\sqrt{\epsilon}}$. Assuming each block uses a scalar regularization, i.e. $c_{i} = \frac{1}{\alpha_{i}}$ where $\alpha_{i} > 0$, we have the following lemma

\begin{lemma} \label{lem:regeigcond}
Let Assumptions~\ref{asm:Qsymmetric} and~\ref{asm:Qbdd} hold for the matrix $Q$ with respect to the partitioning vector $p = [n_{1},n_{2},\dots,n_{N}]^{T}$. Let $A=\left[A^{[i]}_{j}\right]_{p}$, with $A^{[i]}_{i} = \alpha_{i}I$ for every $i \in [N]$ and $A^{[i]}_{j} = 0$ when $j \neq i$. If we have
$\alpha_{i} \leq \frac{\sqrt{\epsilon}}{1-\sqrt{\epsilon}}\delta_{i}(Q) \textnormal{ for all } i \in [N]$,
    then $\lambda_{min}\left(A^{-1}Q\right) \geq \frac{1-\sqrt{\epsilon}}{\sqrt{\epsilon}}$.
\end{lemma}

\textit{Proof:} Use Equation~\eqref{eqn:CB>l} and substitute $C=A^{-1}$, $B=Q$, and $l = \frac{1-\sqrt{\epsilon}}{\sqrt{\epsilon}}$. $\hfill\blacksquare$

We have shown this eigenvalue condition can be satisfied according to the conditions in Problem 3, i.e. $A^{[i]}$ is chosen using only knowledge of $Q^{[i]}$ and $\epsilon$. The following subsection will show this condition is sufficient to satisfy the error bound in Problem 3.

\subsection{Error Bound Satisfaction}

Proof of error bound satisfaction will be done using the following lemma.

\begin{lemma} \label{lem:AQtoError}
Let $f(x) = \frac{1}{2}x^{T}Qx + r^{T}x$, where $Q=Q^{T} \succ 0$, $Q \in \mathbb{R}^{n \times n}$, and $r, x \in \mathbb{R}^{n}$. Let $\hat{x} = \textnormal{arg}\min_{x \in \mathbb{R}^{n}} f(x)$ and $\hat{x}_{A} = \textnormal{arg}\min_{x \in \mathbb{R}^{n}} f(x)+\frac{1}{2}x^{T}Ax$, where $A \succ 0$ and diagonal. Additionally, let $\epsilon \in [0,1]$. If 
\begin{equation}
   \frac{1-\sqrt{\epsilon}}{\sqrt{\epsilon}} \leq \lambda_{min}(A^{-1}Q)\textnormal{, then }
    \frac{\left| f(\hat{x}) - f(\hat{x}_{A})\right|}{\left| f(\hat{x})\right|} \leq \epsilon.
\end{equation}

\end{lemma}

\textit{Proof:} Proof in Appendix~\ref{app:AQtoError}. $\hfill\blacksquare$

With these lemmas, we now present the following theorem.

\begin{theorem} \label{the:errorbound}
Let Assumptions~\ref{asm:Qsymmetric} and~\ref{asm:Qbdd} hold for the matrix $Q$ with respect to the partitioning vector $p = [n_{1},n_{2},\dots,n_{N}]^{T}$. Let $A=\left[A^{[i]}_{j}\right]_{p}$, with $A^{[i]}_{i} = \alpha_{i}I$ for every $i \in [N]$ and $A^{[i]}_{j} = 0$ when $j \neq i$. Let $f(x) = \frac{1}{2}x^{T}Qx + r^{T}x$, where $r, x \in \mathbb{R}^{n}$. Let $\hat{x} = \textnormal{arg}\min_{x \in \mathbb{R}^{n}} f(x) = -Q^{-1}r$ and $\hat{x}_{A} = \textnormal{arg}\min_{x \in \mathbb{R}^{n}}f(x)+\frac{1}{2}x^{T}Ax = -P^{-1}r$, where $P=Q+A$. Additionally, let $\epsilon \in [0,1]$. If
\begin{equation} \label{eqn:ealphabound}
\alpha_{i} \leq \frac{\sqrt{\epsilon}}{1-\sqrt{\epsilon}}\delta_{i}(Q) \textnormal{ for all } i \in [N],
\end{equation}
then,
\begin{align}
    \frac{\left| f(\hat{x}) - f(\hat{x}_{A})\right|}{\left| f(\hat{x})\right|} & \leq \epsilon
\end{align}
\end{theorem}

\textit{Proof:} Lemma~\ref{lem:regeigcond} shows that the regularization selection rules presented above, along with Assumption~\ref{asm:Qbdd}, imply that $\frac{1-\sqrt{\epsilon}}{\sqrt{\epsilon}} \leq \lambda_{min}(A^{-1}Q)$. Lemma~\ref{lem:AQtoError} shows that $\frac{1-\sqrt{\epsilon}}{\sqrt{\epsilon}} \leq \lambda_{min}(A^{-1}Q)$ implies that $\frac{\left| f(\hat{x}) - f(\hat{x}_{A})\right|}{\left| f(\hat{x})\right|} \leq \epsilon$. $\hfill\blacksquare$

Additionally, we can derive a similar bound for relative error in the solution itself. Defining this error as $\frac{\|\hat{x}-\hat{x}_{A}\|_{2,p}}{\|\hat{x}\|_{2,p}}$ and using Equation~\eqref{eqn:x-xaform} we see
\begin{align}
    & \frac{\|\hat{x}-\hat{x}_{A}\|_{2,p}}{\|\hat{x}\|_{2,p}} = \frac{\|(I+A^{-1}Q)^{-1}Q^{-1}r\|_{2,p}}{\|Q^{-1}r\|_{2,p}} \\
    & \leq \frac{\|(I+A^{-1}Q)^{-1}\|_{2,p}\|Q^{-1}r\|_{2,p}}{\|Q^{-1}r\|_{2,p}} = \|(I+A^{-1}Q)^{-1}\|_{2,p} \\
    & \leq \frac{1}{\min_{i \in [N]}\left[1+\alpha^{-1}_{i}\delta_{i}(Q)\right]}.
\end{align}
If we wish for agents to select regularizations such that the above error is less than a given constant $\eta$, we see this is accomplished if
\begin{align}
    \frac{1}{\eta} & \leq \min_{i \in [N]}1+\alpha^{-1}_{i}\delta_{i}(Q) \\
    \alpha_{i} & \leq \frac{\eta}{1-\eta}\delta_{i}(Q) \textnormal{ for all } i \in [N].
\end{align}
This rule has the same structure as the one in Theorem~\ref{the:errorbound}, with the only difference being there is no square root taken of $\eta$. 

Note that throughout this section it was assumed that $A$ is invertible, which is true if $\alpha_{i} > 0$ for all $i \in [N]$. However in scenarios where there is no need for a particular agent to regularize, e.g. $q_{i} < q^{*}$, that agent can choose $\alpha_{i} = 0$ for all practical applications. This is because all of the above analysis holds if $\alpha_{i}$ is chosen to be a small positive value, which can be set arbitrarily close to zero.

\subsection{Trade-Off Analysis}
There is an inherent trade-off between the speed at which we reach a solution and the quality of that solution. Theorem~\ref{thm:qreg} provides a lower bound on $\alpha_{i}$ that allows us to converge at any speed we wish, while Theorem~\ref{the:errorbound} provides an upper bound on $\alpha_{i}$ that allows us to bound the cost error between the solution we find and the optimal solution. However, in general, there is no reason to expect these two bounds to be compatible in the sense that $\alpha_{i}$ can be chosen such that both are satisfied for all $i \in [N]$. Therefore, when implemented, it is likely that the network operator will be able to choose whether speed or accuracy is more critical for the specific problem. If speed is mission-critical, then agents may select the smallest regularizations required to match that speed, and if accuracy is mission-critical, agents may select the largest regularizations that obey the specific error bound.

\section{Simulation} \label{sec:simulation}
\begin{figure}[!tp]
\centering
\includegraphics[draft = false,width=3.6in]{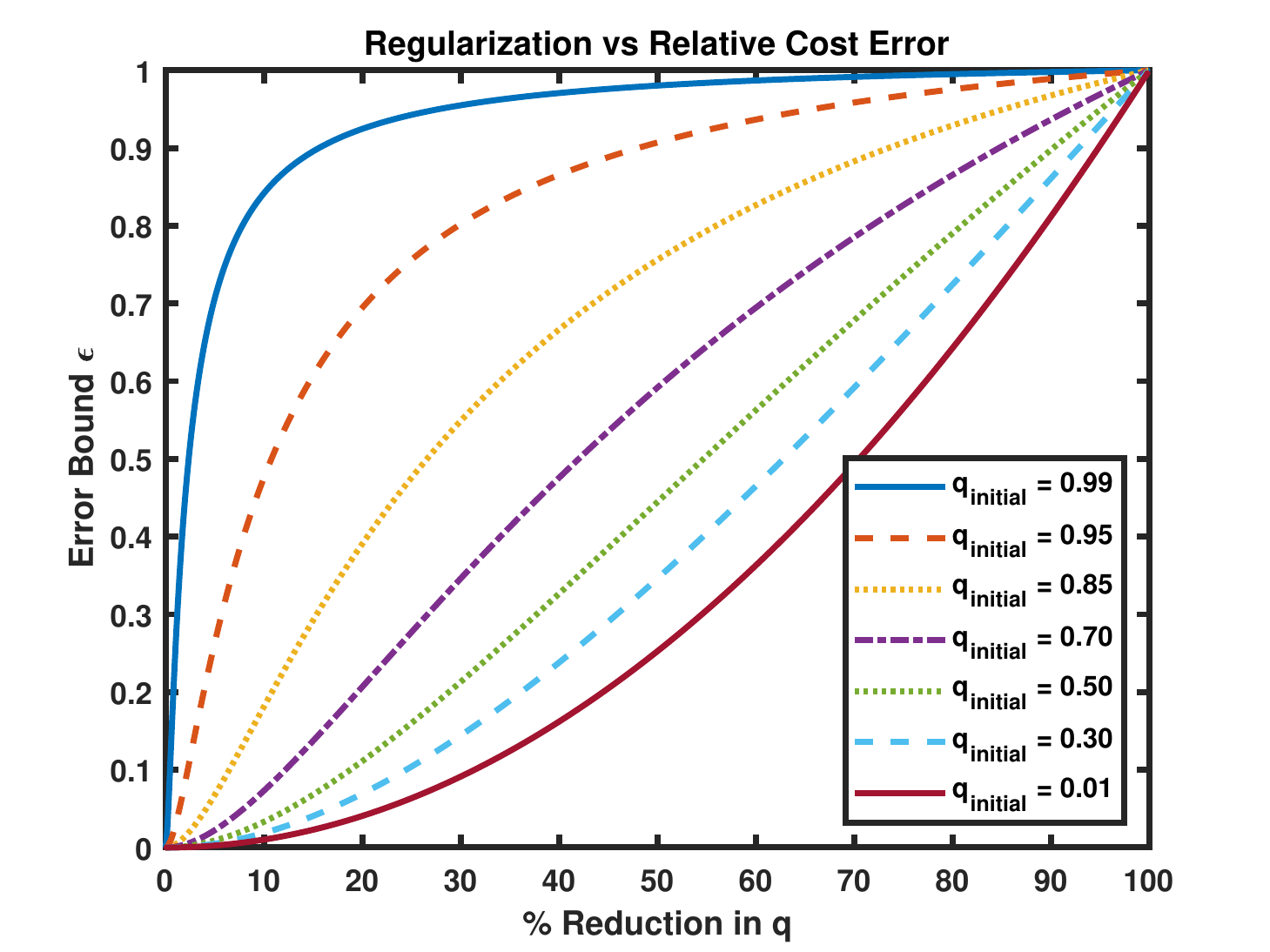}
\caption{The percent reduction in $q$ due to regularization plotted vs the relative cost error bound that regularization induces, with different lines plotting this relationship for QPs with different initial values for $q$.
}
\label{fig:qvse} 
\end{figure}

To visualize the trade-off between speed and error when regularizing, we generate seven QPs, each with 100 diagonally dominant blocks. The QPs are generated to have initial convergence parameters of $q_{initial} =$ 0.99, 0.95, 0.85, 0.70, 0.50, 0.30, and 0.01. For each QP, $A$ is independently chosen according to Theorem~\ref{thm:qreg} such that $q$ is reduced by percentages ranging from 0\% to 100\%, and this percentage reduction is plotted against the corresponding error bound given by Theorem~\ref{the:errorbound} in Figure~\ref{fig:qvse}. For example, the data for the QP with $q_{initial} = 0.85$ is plotted by the yellow dotted line in Figure~\ref{fig:qvse}, and one can see that if this QP is regularized to reduce $q$ by 10\% (i.e., a reduction from 0.85 to 0.765), the relative error in cost can be upper bounded by approximately $\epsilon = 18\%$.

There are two main takeaways from Figure~\ref{fig:qvse}. The first is that, as expected, larger regularizations result in a larger relative error bound, which is upper bounded by 1. This is because $q \rightarrow 0$ as $A \rightarrow \infty$, $f(\hat{x}_{A}) \rightarrow 0$ as $A \rightarrow 0$, and $\epsilon \rightarrow 1$ as $f(\hat{x}_{A}) \rightarrow 0$. The second is that the larger $q_{initial}$ is, the more sensitive the error bound for the QP is to regularizing. That is, if $q_{initial}$ is thought of as a condition number, then ``poorly conditioned" QPs will have larger errors due to regularizing.

\begin{figure}[!tp]
\centering
\includegraphics[draft=false,width=3.6in]{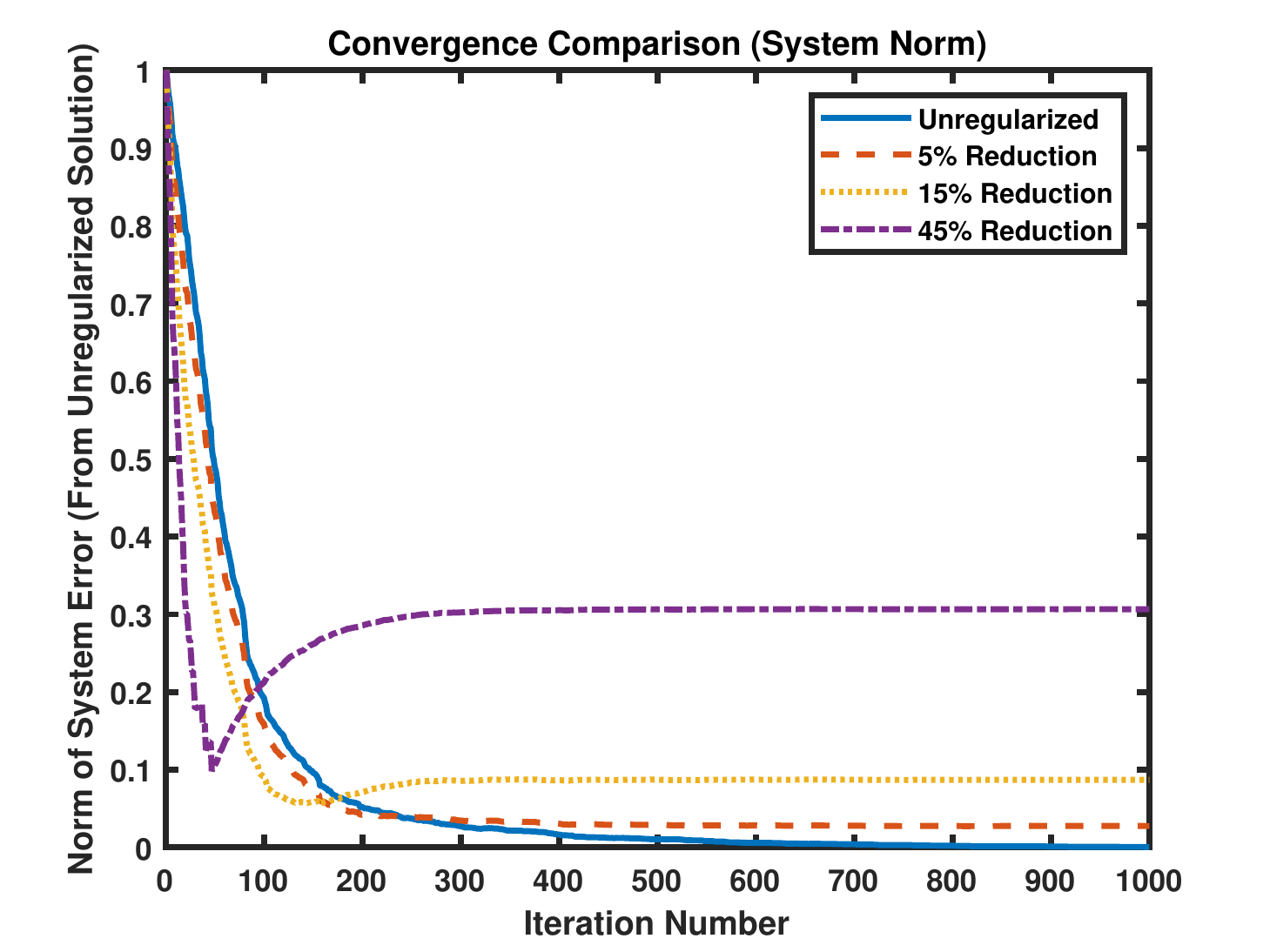}
\caption{Network error convergence of Algorithm 1 when unregularized vs regularizing such that $q$ is reduced by 5\%, 15\%, and 45\%.
}
\label{fig:comparison} 
\end{figure}

A second simulation was run to demonstrate the convergence properties due to regularizing. One QP was generated with 100 blocks and $q_{initial} = 0.85$. Three different regularization matrices were chosen according to Theorem~\ref{thm:qreg}, called $A_{5}$, $A_{15}$, and $A_{45}$, such that $q$ is reduced by 5\%, 15\%, and 45\%, respectively. The blocks are then distributed among 100 agents, who have a 10\% chance of computing an update and a 1\% chance of transmitting a state to each other agent at each timestep. Four simulations were run, one solving the unregularized QP, and three others using each regularization matrix. The 2-norm of the system error to the unregularized solution, $\|x(k)-\hat{x}\|_{2}$, is plotted for each simulation against iteration number in Figure~\ref{fig:comparison}.

As expected, only the unregularized case converges to the unregularized solution, while the other cases converge to other solutions whose distances to the unregularized solution grow with larger regularizations. However, the cases with larger regularizations initially converge to $\hat{x}$ faster, up to a point. That is, larger regularizations mean the system will initially move toward $\hat{x}$ faster, but will reach the turn-off point, where the system error grows again, earlier and further away from $\hat{x}$. This behavior suggests a vanishing regularization scheme, where $A$ shrinks to zero with time, may lead to accelerated convergence to the exact solution $\hat{x}$. Note also that convergence even in the unregularized case is non-monotone, and at times the norm of the system error may even grow due to the asynchronous nature of of communications, but Theorem~\ref{thm:alg1works} guarantees these growths are bounded and error will converge to zero.

\section{Conclusions} \label{sec:conclusions}
We have developed a distributed quadratic programming framework that converges under totally asynchronous conditions. This framework allows agents to select stepsizes and regularizations independently of one another, using only knowledge of their block of the QP, that guarantee a specified global convergence rate and cost error bound. Future work will apply these developments to quadratic programs with functional constraints.

\bibliographystyle{IEEEtran}
\bibliography{Biblio}

\begin{IEEEbiography}[{\includegraphics[width=1in,height=1.25in,clip,keepaspectratio,draft=false]{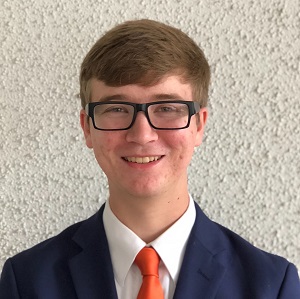}}]{Matthew Ubl} is a PhD student at the University of Florida, where he is an Institute for Networked Autonomous Systems Fellow and a recipient of the Graduate Student Preeminence Award. He received his Bachelor's degree in Aerospace Engineering from the University of Central Florida in 2018. His current research interests are in asynchronous multi-agent coordination, with particular focus upon multi-agent optimization with impaired and unreliable communications.
\end{IEEEbiography}
\begin{IEEEbiography}[{\includegraphics[width=1in,height=1.25in,clip,keepaspectratio,draft=false]{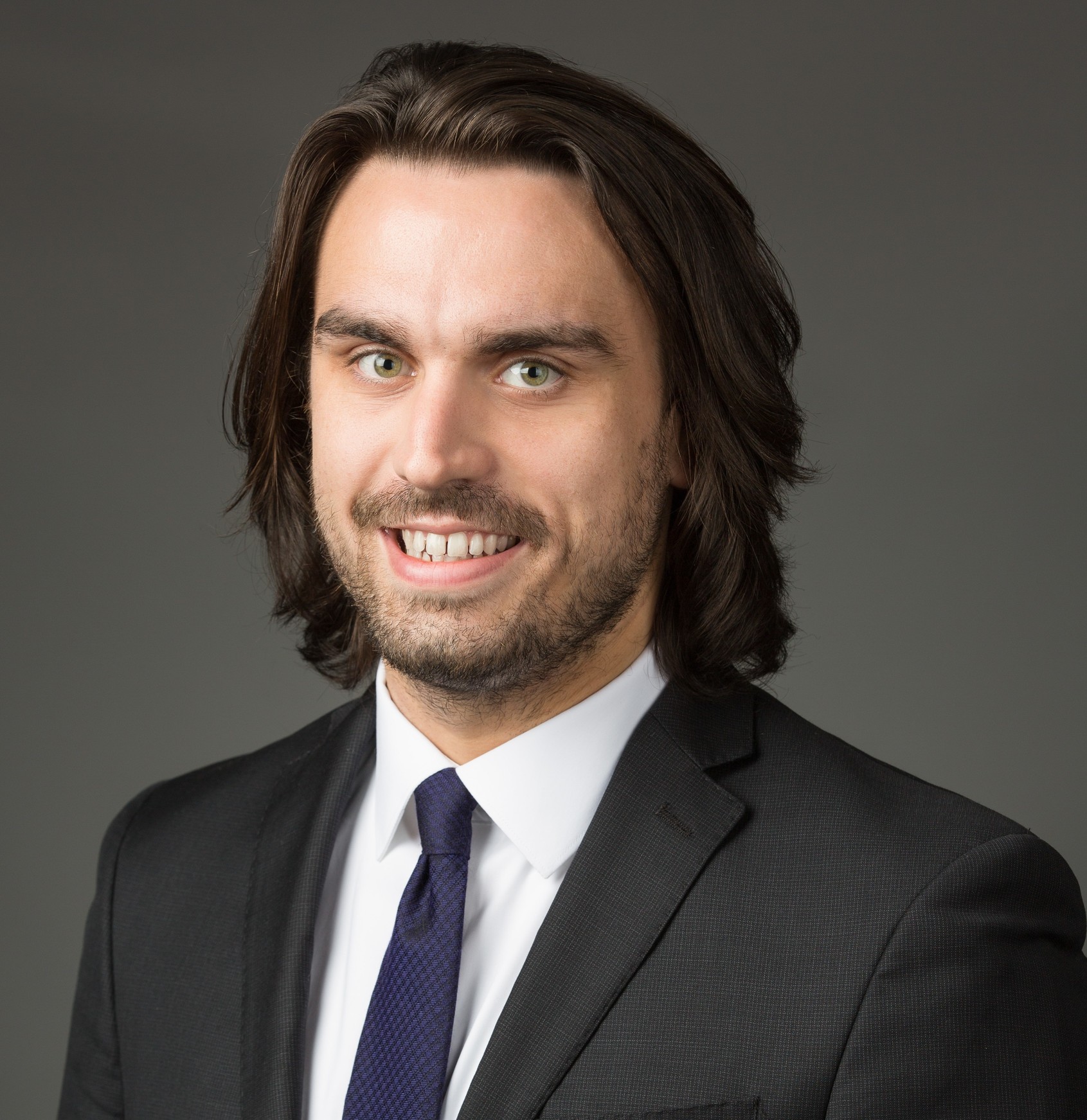}}] {Matthew Hale} is an Assistant Professor of Mechanical and Aerospace Engineering at the University of Florida. He received his BSE in Electrical Engineering summa cum laude from the University of Pennsylvania in 2012, and his MS and PhD in Electrical and Computer Engineering from the Georgia Institute of Technology in 2015 and 2017, respectively. He directs the Control, Optimization, and Robotics Engineering (CORE) Lab at the University of Florida, and his research interests include multi-agent systems, mobile robotics, privacy in control, distributed optimization, and graph theory. He was the Teacher of the Year in the UF Department of Mechanical and Aerospace Engineering for the 2018-2019 school year, and he received an NSF CAREER Award in 2020 for his work on privacy in control systems.   
\end{IEEEbiography}

\appendices
\section{} \label{app:normsumbound}
Proof of Lemma~\ref{lem:normsumbound}: By definition of the maximum norm,
\begin{align}
    \|B\|_{2,p} & = \sup_{\|x\|_{2,p}=1}\|Bx\|_{2,p} 
                     = \sup_{\|x\|_{2,p}=1}\max_{i \in [N]}\left\|B^{[i]}x\right\|_{2}.
\end{align}

Since $B^{[i]} = \left[B^{[i]}_1 \textnormal{ }B^{[i]}_2 \textnormal{ ... } B^{[i]}_N\right]$, we can now write ${B^{[i]}x = B^{[i]}_{1}x^{[1]} \textnormal{ }+\textnormal{ }B^{[i]}_{2}x^{[2]}+ \textnormal{ ... } +B^{[i]}_{N}x^{[N]}}$. Therefore,
\begin{align}
    \|B\|_{2,p} & = \sup_{\|x\|_{2,p}=1}\max_{i \in [N]}\left\|B^{[i]}_{1}x^{[1]} \textnormal{ }+ \textnormal{ ... } +B^{[i]}_{N}x^{[N]}\right\|_{2}.
\end{align}

By the triangle inequality, we have
\begin{align}
    \|B\|_{2,p} & \leq \sup_{\|x\|_{2,p}=1}\max_{i \in [N]}\sum^{N}_{j=1}\left\|B^{[i]}_{j}x^{[j]}\right\|_{2}.
\end{align}

The condition $\|x\|_{2,p} = 1$ implies $\left\|x^{[i]}\right\|_{2} \leq 1$ for all $i \in [N]$. Therefore, for each element in the sum above, we can write ${\left\|B^{[i]}_{j}x^{[j]}\right\|_{2} \leq \sup_{\left\|x^{[j]}\right\|_{2}=1}\left\|B^{[i]}_{j}x^{[j]}\right\|_{2} =  \left\|B^{[i]}_{j}\right\|_{2}}$. Substituting this above completes the proof. $\hfill\blacksquare$

\section{} \label{app:ddnec}
Proof of Lemma~\ref{lem:ddnec}:
Because $Q^{[i]}_{i}=Q^{[i]^{T}}_{i} \succ 0$, we see that 
\begin{align}
    & \left\|I -\gamma_{i}Q^{[i]}_{i}\right\|_{2} \\
    & = \textnormal{max}\left\{ \left|\lambda_{min}\left(I-\gamma_{i} Q^{[i]}_{i}\right)\right|,\left|\lambda_{max}\left(I-\gamma_{i} Q^{[i]}_{i}\right)\right|\right\} \\
    & = \textnormal{max}\left\{ \left|1-\gamma_{i} \lambda_{min}\left(Q^{[i]}_i\right)\right|,\left|1-\gamma_{i} \lambda_{max}\left(Q^{[i]}_i\right)\right| \right\},
\end{align}
which allows us to write 
\begin{align}
    \left\|I-\gamma_{i}Q^{[i]}_{i}\right\|_{2} + \gamma_{i}\sum^{N}_{\substack{j=1 \\ j\neq i}}\left\|Q^{[i]}_{j}\right\|_{2} < 1
\end{align}
if and only if both
\begin{align}
    \left|1-\gamma_{i} \lambda_{min}\left(Q^{[i]}_i\right)\right| < 1 - \gamma_{i}\sum^{N}_{\substack{j=1 \\ j\neq i}}\left\|Q^{[i]}_{j}\right\|_{2}
    \end{align}
and
\begin{align}
    \left|1-\gamma_{i} \lambda_{max}\left(Q^{[i]}_i\right)\right| < 1 - \gamma_{i}\sum^{N}_{\substack{j=1 \\ j\neq i}}\left\|Q^{[i]}_{j}\right\|_{2}.
\end{align}

The first inequality will be true if and only if both
\begin{equation} \label{eqn:BDDconditionsatisfied}
    \lambda_{min}(Q^{[i]}_i)> \sum^{N}_{\substack{j=1 \\ j\neq i}}\left\|Q^{[i]}_{j}\right\|_{2}
    \end{equation}
    and 
    \begin{equation}
\gamma_{i} < \frac{2}{\lambda_{min}(Q^{[i]}_i) + \sum^{N}_{\substack{j=1 \\ j\neq i}}\left\|Q^{[i]}_{j}\right\|_{2}},
\end{equation}
and the second will be true if and only if both
\begin{equation}
    \lambda_{max}(Q^{[i]}_i)> \sum^{N}_{\substack{j=1 \\ j\neq i}}\left\|Q^{[i]}_{j}\right\|_{2}
    \end{equation}
    and 
    \begin{equation} \label{SSconditionsatisfied}
\gamma_{i} < \frac{2}{\lambda_{max}(Q^{[i]}_i) + \sum^{N}_{\substack{j=1 \\ j\neq i}}\left\|Q^{[i]}_{j}\right\|_{2}}.
\end{equation}

Taking the most restrictive of these conditions, we can write
\begin{align}
    \left\|I-\gamma_{i}Q^{[i]}_{i}\right\|_{2} + \gamma_{i}\sum^{N}_{\substack{j=1 \\ j\neq i}}\left\|Q^{[i]}_{j}\right\|_{2} < 1
    \end{align}
    if and only if Equations~\eqref{eqn:BDDconditionsatisfied} and~\eqref{SSconditionsatisfied} hold. $\hfill\blacksquare$

\section{} \label{app:setswork}   
Proof of Theorem~\ref{thm:setswork}: For Assumption~\ref{asm:setsexist}.1, by definition we have
\begin{equation}
X(s+1)=\left\{ y\in X:\|y-\hat{x}\|_{2,p}\leq q^{s+1}D_{o}\right\} .
\end{equation}
Since $q\in(0,1)$, we have $q^{s+1}<q^{s}$, which results in~${\|y-\hat{x}\|_{2,p}\leq q^{s+1}D_{o}<q^{s}D_{o}}$.
Then $y\in X(s+1)$ implies $y\in X(s)$ and $X(s+1)\subset X(s)\subset X$,
as desired.

For Assumption~\ref{asm:setsexist}.2 we find
\begin{equation}
\lim_{s\rightarrow\infty}X(s)=\lim_{s\rightarrow\infty}\left\{ y\in X:\|y-\hat{x}\|_{2,p}\leq q^{s}D_{o}\right\} =\left\{ \hat{x}\right\}.
\end{equation}
The structure of the weighted block-maximum
norm then allows us to see that $\|y-\hat{x}\|_{2,p}\leq q^{s}D_{o}$
if and only if $\|y^{[i]}-\hat{x}^{[i]}\|_{2}\leq q^{s}D_{o}$
for all $i\in[N].$ It then follows that
\begin{equation}
X_{i}(s)=\left\{ y^{[i]}\in X_{i}:\|y^{[i]}-\hat{x}^{[i]}\|_{2}\leq q^{s}D_{o}\right\} ,
\end{equation}
which gives $X(s)=X_{1}(s)\times...\times X_{N}(s)$, thus satisfying
Assumption~\ref{asm:setsexist}.3.

We then see that, for $y\in X(s)$,
\begin{align}
\left\|\theta_{i}(y)-\hat{x}^{[i]}\right\|_{2}= \Big\|\Pi_{X_{i}}&\left[y^{[i]}- \gamma_{i}\left(Q^{[i]}y+r^{[i]}\right)\right] \\
& -\Pi_{X_{i}}\left[\hat{x}^{[i]}-\gamma_{i}\left(Q^{[i]}\hat{x}+r^{[i]}\right)\right]\Big\|_{2},
\end{align}
which follows from the definition of $\theta_{i}(y)$ and the fact
that $\hat{x}^{[i]}=\Pi_{X_{i}}\left[\theta_{i}(\hat{x})\right]$. Using the non-expansive property of the projection operator $\Pi_{X_{i}}\left[\cdot\right]$, we find
\begin{align*}
\left\|\theta_{i}(y)-\hat{x}^{[i]}\right\|_{2} & \leq \left\|y^{[i]}-\hat{x}^{[i]}-\gamma_{i}Q^{[i]}\left(y -\hat{x}\right)\right\|_{2}\\
& = \left\|\left(I^{[i]}-\gamma_{i}Q^{[i]}\right)\left(y -\hat{x}\right)\right\|_{2}\\
& \leq \max_{i\in[N]} \left\|\left(I^{[i]}-\gamma_{i}Q^{[i]}\right)\left(y -\hat{x}\right)\right\|_{2}\\
& = \left\|\left(I-\Gamma Q\right)\left(y -\hat{x}\right)\right\|_{2,p} \\
& \leq \|I-\Gamma Q\|_{2,p} \|y -\hat{x}\|_{2,p},
\end{align*}
which follows from our definition of the block-maximum norm. From Lemmas~\ref{lem:normsumbound} and~\ref{lem:ddnec} we know ${\|I-\Gamma Q\|_{2,p} \leq q < 1}$,
and using the hypothesis that ${y\in X(s)}$, we find
\begin{align*}
\left\|\theta_{i}(y)-\hat{x}^{[i]}\right\|_{2} & \leq q \|y -\hat{x}\|_{2,p} \leq q^{s+1} D_{o},
\end{align*}
which shows $\theta_{i}(y)\in X_{i}(s+1)$ and
Assumption~\ref{asm:setsexist}.4 is satisfied. $\hfill\blacksquare$

\section{} \label{app:alg1works}
Proof of Theorem~\ref{thm:alg1works}: Theorem~\ref{thm:setswork} shows the construction of the sets $\left\{ X(s)\right\} _{s\in\mathbb{N}}$
satisfies Assumption~\ref{asm:setsexist}, and from ~\cite{bertsekas1989convergence} and ~\cite{bertsekas1989parallel} we see this implies asymptotic
convergence of Algorithm 1 for all $i\in[N]$. The total
asynchrony required by Problem 1 is incorporated by not requiring
delay bounds, and agents do not require any coordination in selecting
stepsizes because the bound on $\gamma_{i}$ depends only upon $Q^{[i]}$, which means that all of the criteria of Problem 1 are satisfied.
$\hfill\blacksquare$

\section{} \label{app:convergencerate}
Proof of Theorem~\ref{convergencerate}:
From the definition of $D_{o}$, for all $i \in [N]$ we have $x_{i}(0) \in X(0)$. If agent $i$ computes a state update, then $\theta_{i}(x_{i}(0)) \in X_{i}(1)$ and after one cycle is completed, say at time $k$, we have $x_{i}(k) \in X(1)$ for all $i$. Iterating this process, after $c(k)$ cycles have been completed by some time $k$, $x_{i}(k) \in X(c(k))$. The result follows by expanding the definition of $X\left(c(k)\right)$. $\hfill\blacksquare$

\section{} \label{app:optimalstep}
Proof of Theorem~\ref{thm:optimalstep}:
If $\gamma_{i} \leq \frac{2}{\lambda_{max}\left(Q^{[i]}_i\right)+\lambda_{min}\left(Q^{[i]}_i\right)},$ then
\begin{align*}
    q_{i} & = 1-\gamma_{i} \left(\lambda_{min}\left(Q^{[i]}_i\right) -\sum^{N}_{\substack{j=1 \\ j\neq i}}\left\|Q^{[i]}_{j}\right\|_{2}\right),
\end{align*}
and if $\gamma_{i} \geq \frac{2}{\lambda_{max}\left(Q^{[i]}_i\right)+\lambda_{min}\left(Q^{[i]}_i\right)},$ then
\begin{align*}
    q_{i} & = - 1 + \gamma_{i} \left(\lambda_{max}\left(Q^{[i]}_i\right) +\sum^{N}_{\substack{j=1 \\ j\neq i}}\left\|Q^{[i]}_{j}\right\|_{2}\right).
\end{align*}

That is, when $\gamma_{i} \leq \frac{2}{\lambda_{max}\left(Q^{[i]}_i\right)+\lambda_{min}\left(Q^{[i]}_i\right)},$ the relationship between $q_{i}$ and $\gamma_{i}$ is a line with negative slope, and when $\gamma_{i} \geq \frac{2}{\lambda_{max}\left(Q^{[i]}_i\right)+\lambda_{min}\left(Q^{[i]}_i\right)}$ the relationship is a line with positive slope. Then $q_{i}$ is minimized at the point where the slope changes sign, which occurs when
\begin{equation} \label{eqn:optstep}
    \gamma_{i} = \frac{2}{\lambda_{max}\left(Q^{[i]}_i\right)+\lambda_{min}\left(Q^{[i]}_i\right)}.
\end{equation}

If every $q_{i}$ has been minimized, then by definition $q$ has been minimized. $\hfill\blacksquare$

\section{} \label{app:AQtoError}
Proof of Lemma~\ref{lem:AQtoError}:
To facilitate this proof, we first present the following facts to which we will repeatedly refer:

\textit{Fact 1:} If $B$ is a square matrix such that $0 < \lambda_{min}(B) \leq \lambda_{max}(B)$, then $\lambda_{max}(B^{-1}) = \lambda_{min}^{-1}(B)$.

\textit{Fact 2:} If $B$ is a square matrix such that $0 < \lambda_{min}(B) \leq \lambda_{max}(B)$, then $\lambda_{min}(B^{2}) = \lambda^{2}_{min}(B)$.

\textit{Fact 3:} If $B$ is a square matrix, then $-\lambda_{max}(B) = \lambda_{min}(-B)$.

\textit{Fact 4:} If $B$ is a square matrix and $C$ is an invertible matrix of the same dimension, then $\lambda_{i}(C^{-1}BC) = \lambda_{i}(B)$ for all $i$.

\textit{Fact 5:} If $B=B^{T} \preceq 0$ and $C$ is an invertible matrix of the same dimension, then $\lambda_{i}(C^{T}BC) \leq 0$ for all $i$.

Facts 1-3 can be easily shown, Fact 4 simply states eigenvalues are invariant under a similarity transform, and Fact 5 is a corollary of Sylvester's Law of Inertia ~\cite[Fact 5.8.17]{bernstein2009matrix}.

Bearing these facts in mind, we first rearrange the condition in the lemma statement to find
\begin{align}
    \frac{1}{\sqrt{\epsilon}} -1 & \leq \lambda_{min}(A^{-1}Q) \\
    \frac{1}{\sqrt{\epsilon}} & \leq 1+ \lambda_{min}(A^{-1}Q) = \lambda_{min}(I+A^{-1}Q) \\
     & = \lambda_{min}(A^{-1}(A+Q)) = \lambda_{min}(A^{-1}P)\\
    \lambda^{-1}_{min}(A^{-1}P) & \leq \sqrt{\epsilon}.
\end{align}
From Fact 1, it follows that $\lambda_{max} (P^{-1} A) \leq \sqrt{\epsilon}$ and
$\lambda^{2}_{max} (P^{-1} A) \leq \epsilon$. From Fact 2, $\lambda_{max} ((P^{-1} A)^2) \leq \epsilon$, which implies $- \epsilon \leq -\lambda_{max} ((P^{-1} A)^2)$. From Fact 3,
\begin{align}
    - \epsilon & \leq \lambda_{min} (-(P^{-1} A)^2) \\
    1- \epsilon & \leq 1+ \lambda_{min} (-(P^{-1} A)^2) = \lambda_{min} (I-(P^{-1} A)^2) \\
    & = \lambda_{min} ((I+P^{-1} A)(I-P^{-1} A)).
\end{align}

Note that $I-P^{-1}A = P^{-1}(P-A) = P^{-1}Q$, therefore
\begin{align}
    1- \epsilon & \leq \lambda_{min} ((I+P^{-1} A)P^{-1}Q) \\
    1- \epsilon & \leq \lambda_{min} ((P^{-1}+P^{-1} AP^{-1})Q).
\end{align}

Note that $P^{-1}+P^{-1} AP^{-1} = P^{-1}+P^{-1} (P-Q) P^{-1} = 2P^{-1} - P^{-1} QP^{-1}$, therefore $1- \epsilon \leq \lambda_{min} ((2P^{-1} - P^{-1} QP^{-1})Q)$,
which implies
\begin{align}
    0 & \leq -(1-\epsilon ) + \lambda_{min} ((2P^{-1} - P^{-1} QP^{-1})Q) \\
    0 & \leq \lambda_{min} (-(1-\epsilon )I + (2P^{-1} - P^{-1} QP^{-1})Q).
\end{align}

From Fact 3, $0 \leq -\lambda_{max} ((1-\epsilon )I - (2P^{-1} - P^{-1} QP^{-1})Q)$
and $\lambda_{max} ((1-\epsilon )I - (2P^{-1} - P^{-1} QP^{-1})Q) \leq 0$.

From Fact 4, taking $C = Q^{-\frac{1}{2}}$
\begin{align}
    \lambda_{max} ((1-\epsilon )I - Q^{\frac{1}{2}}(2P^{-1} - P^{-1} QP^{-1})Q^{\frac{1}{2}}) & \leq 0.
\end{align}

Note that the matrix above is symmetric. Therefore, from Fact 5, taking $C = Q^{-\frac{1}{2}}$, we have
\begin{align}
    \lambda_{max} ((1-\epsilon )Q^{-1} - 2P^{-1} + P^{-1} QP^{-1}) & \leq 0.
\end{align}

Note that the matrix above is still symmetric. Therefore, we can write 
$(1-\epsilon )Q^{-1} - 2P^{-1} + P^{-1} QP^{-1} \preceq 0$,
which implies $Q^{-1} - 2P^{-1} + P^{-1} QP^{-1}  \preceq \epsilon Q^{-1}$.

This means that for any arbitrary vector $x$ of dimension $n$, 
$x^{T}(Q^{-1} - 2P^{-1} + P^{-1} QP^{-1})x \leq x^{T}(\epsilon Q^{-1})x$,
and
$x^{T}Q^{-1}x - 2x^{T}P^{-1}x + x^{T}P^{-1} QP^{-1}x \leq \epsilon x^{T}Q^{-1}x$.

Assuming $x \neq 0$, $x^{T}Q^{-1}x$ is a positive scalar. Dividing both sides by this term gives
\begin{align}
    \frac{x^{T}Q^{-1}x - 2x^{T}P^{-1}x + x^{T}P^{-1} QP^{-1}x}{x^{T}Q^{-1}x} & \leq \epsilon.
\end{align}

Because this relation is true for any arbitrary vector, we can choose $x=r$ and multiply by $\frac{-\frac{1}{2}}{-\frac{1}{2}}$ to find
\begin{align}
    \frac{-\frac{1}{2}r^{T}Q^{-1}r - (\frac{1}{2}r^{T}P^{-1} QP^{-1}r - r^{T}P^{-1}r)}{-\frac{1}{2}r^{T}Q^{-1}r} & \leq \epsilon,
\end{align}
and substituting returns the desired result. $\hfill\blacksquare$
\end{document}